\input amstex
\documentstyle{amsppt}
\baselineskip12pt
\magnification=1100
\centerline{\bf Inverse Spectral Problems in Rectangular Domains}
\medskip
\centerline{Gregory Eskin and James Ralston, UCLA}
\centerline{Los Angeles, CA 90095-1555, USA}
\centerline{email: eskin\@math.ucla.edu, ralston\@math.ucla.edu}
\bigskip
\noindent {\bf Abstract.} We consider the Schr\"odinger operator 
$-\Delta +q$ in domains of the form
$R=\{x\in {\Bbb R}^n:\ 0\leq x_i\leq a_i,\ i=1,..,n\}$ with either 
Dirichlet or Neumann boundary conditions on the faces 
of $R$, and study the constraints on $q$ imposed by fixing the spectrum 
of $-\Delta +q$ with these boundary conditions. We work in the space of 
potentials, $q$, which become real-analytic on ${\Bbb R}^n$ when they 
are extended evenly across the coordinate planes and then periodically. 
Our results have the corollary that there are no continuous isospectral 
deformations for these operators within that class of potentials.
This work is based on new formulas for the trace of the wave group in 
this setting. In addition to the inverse
spectral results these formulas lead to asymptotic expansions for the 
traces of the wave and heat kernels on rectangular domains.
\medskip
\noindent {\bf \S 1. Introduction}
\medskip

We consider the Schr\"odinger operator, $H=-\Delta +q(x)$, in the 
rectangular domain
$R=\{x\in {\Bbb R}^n:\ 0\leq x_i\leq a_i,\ i=1,..,n\}$
and on each of the $2^n$ faces of $\partial R$, $R\cap\{x_i=0\hbox{ or 
}x_i=a_i\}$ we impose either Dirichlet or Neumann boundary conditions.
Let
$$\mu_1<\mu_2\leq\mu_3\leq \cdots $$
be the spectrum of the self-adjoint operator, also denoted by $H$, that 
is obtained in this way.
We take $R$ and $\{\mu_j\}_{j=1}^\infty$ as given, and study what 
constraints
these data put on $q$.

This work depends strongly on the geometry of rectangular domains, and 
involves many reflections, translations, etc.
Since all of this will be easier to follow once one has seen it in a 
simple case, we have moved the statements and proofs of the general 
results to \S 5. In this Introduction and \S 2-\S 4 we will restrict 
ourselves to the case $n=2$ with Dirichlet boundary conditions, i.e.  we will 
work in the rectangle
$$R=\{(x_1,x_2):0\leq x_1\leq a,\ 0\leq x_2\leq b\},$$
with the boundary condition $u=0$ on $\partial R$.

For our results we need to extend $q$ to an even periodic potential on 
${\Bbb R}^2$:
we extend $q$ to $Q$ on ${\Bbb R}^2$ by defining
$$Q(-x_1,x_2)=q(x_1,x_2)\hbox{ for } x\in R,$$
$$Q(x_1,-x_2)=Q(x_1,x_2)\hbox{ for } |x_1|\leq a,\ 0\leq x_2\leq 
b,\hbox{ and }$$
$$Q(x_1+2ma,x_2+2nb)=Q(x_1,x_2)\hbox{ for } |x_1|\leq a,\ |x_2|\leq 
b,\hbox{ and } (m,n)\in {\Bbb Z}^2.$$ 
Thus $Q$ is periodic with respect to the lattice $L
=\{(2ma,2nb),\ m,n\in { \Bbb Z}\}$. To $L$ we associate the dual 
lattice $$L^*=\{\delta\in {\Bbb R}^2:\delta\cdot d\in {\Bbb Z}\hbox{ for all 
}d\in L\},$$ and expand $Q$ in a Fourier series
$$Q(x)=\sum_{\delta \in L^*}a_\delta e^{2\pi i \delta\cdot x}. $$
Let $S$ denote the set of elements of $S$ which are maximal in the 
sense that $\{\delta\cdot d,\ d\in L\}=\Bbb Z$. For our lattice 
$S=\{(m/(2a),n/(2b)):\hbox{ m and n are relatively prime }\}$. Then we can 
decompose $Q$
$$Q(x)={1\over 2}\sum_{\delta \in S}(\sum_{k=-\infty}^\infty 
a_{k\delta}e^{2\pi i k\delta\cdot x})={1\over 2}\sum_{\delta \in S} 
Q_\delta(\delta\cdot x),$$
where
$$Q_\delta(s)=\sum_{k=-\infty}^\infty a_{k\delta}e^{2\pi i ks}.$$
Note that this decomposition holds if, but only if, 
$a_{(0,0)}=\int_Rqdx=0$. However, in all dimensions (see \S 6) one can recover $\int_Rqdx$ 
from the asymptotics
of the trace $\sum \exp(-\mu_jt)$ as $t\to 0$. So given isospectral 
potentials, subtracting this constant from each of them, one can replace 
them
by isospectral potentials satisfying $\int_Rqdx=0$.
Note also that $Q_{-\delta}(s)=Q_\delta (-s)$ which gives rise to the 
factor $1/2$ in these formulas.
As in [ERT1] we call $Q_\delta$ a \lq\lq directional potential". Since 
$Q(x)=Q(-x)$, we have $Q_\delta(s)=Q_\delta(-s)$. The directional 
potentials corresponding to coordinate directions, i.e.
$\delta_1 =( 1/(2a),0)$ and $\delta_2 =(0,1/(2b))$, only depend on the 
coordinates $x_1$ and $x_2$ respectively, and we denote them by 
$q_1(x_1)$ and $q_2(x_2)$. The main result of the first part of this article 
is:
\medskip
\noindent {\bf Theorem 1.1.} Assume that  $-\Delta +q$ and $-\Delta 
+\tilde q$ have the same Dirichlet spectrum on $R$. If $a^2/b^2$ is 
irrational and the extensions $Q$ and $\tilde Q$ of
$q$ and $\tilde q$ described above are real-analytic on  ${\Bbb R}^2$, 
then
\medskip
\noindent a) for any $\delta\in S$ with no zero components the 
operators $-|\delta|^2d^2/ds^2 +
Q_\delta(s)$ and 

\noindent $-|\delta|^2d^2/ds^2 +
\tilde Q_\delta(s)$ with periodic boundary conditions on $[0,1]$ have 
the same spectrum. 
\medskip
\noindent b) the pairs of operators, $-d^2/dx_1^2 +q_1(x_1)$ and 
$-d^2/dx_1^2 +\tilde q_1(x_1)$, and, $-d^2/dx_2^2 +q_2(x_2)$ and $-d^2/dx_2^2 
+\tilde q_2(x_2)$, have the same spectrum for Dirichlet 
boundary conditions on $[0,a]$ and $[0.b]$ respectively. 
\medskip 
Part a) of Theorem 1.1 was a surprise. In the case of periodic boundary 
conditions it was natural that the directional potentials would be 
isospectral for periodic boundary
conditions, but we did not anticipate that this would be true for 
Dirichlet conditions. Another unexpected result in comparison to the 
periodic case is the rigidity of the directional potentials $Q_\delta$ for 
$\delta$ in the \lq\lq oblique directions" in a). 
Since $Q_\delta$ is even, the isospectral set for $-|\delta|^2d^2/ds^2 
+
Q_\delta(s)$ with periodic boundary conditions on $[0,1]$ is either 
finite or a Cantor set, depending on whether a finite or an infinite 
number of \lq\lq gaps" in the spectrum are open (see [GT]).

While part b) is the most  one can say for pairs of isospectral 
potentials of the form $q(x)=q_1(x_1) +q_2(x_2)$, i.e. for potentials of that 
form b) implies that $q$ and $\tilde q$ are isospectral on $R$, this is 
not the end of the story. Since for our class of potentials, the 
reduced potentials $q_1(x_1)$ and $q_2(x_2)$ also satisfy the one-dimensional
version of the symmetry conditions, it turns their Dirichlet spectra on 
$[0,a]$ and $[0,b]$ determine their periodic spectra on $[-a,a]$ and 
$[-b,b]$ respectively. Thus, like the directional potentials in the 
oblique directions, the set of isospectral potentials in the coordinate 
directions is discrete. In other words {\it within our class of
potentials $q$} the set of isospectral potentials for the Dirichlet 
condition on a rectangle is discrete. 
\medskip
\noindent {\bf Theorem 1.2.} If $a^2/b^2$ is irrational, then any
continuous curve $q_t$ of isospectral potentials such that the 
extensions $Q_t$ are real analytic for all $t$ must be constant,
$q_t=q_0$ for all $t$.
\medskip

The proof of Theorem 1.1 is given in \S 2 - \S 4. It is an extension of 
the proof of the analogous 
result, [ERT1, Theorems 6.1 and 6.2], for periodic boundary conditions. 
The new ingredients in the proof are in the analysis of the trace of 
the fundamental solution of the wave equation in $R$ with the potential 
$q$. The Dirichlet trace is significantly different from the trace for 
periodic
boundary conditions, but the new terms in the Dirichlet trace \lq\lq 
telescope" in a way which
simplifies their contribution to the singularities of the trace. The 
resulting formulas reveal close relations between the Dirichlet and 
periodic traces, and they also make it possible to compute the singularities 
in the trace in a way that identifies contributions with the underlying
geometry. For this we found it simpler to use the classical Hadamard 
construction of parametrices for the wave equation instead of Fourier 
integral operators. Given the lemma that isospectral potentials in the 
coordinate directions form a discrete set, Theorem 1.2 is an immediate 
corollary of Theorem 1.1. The lemma is proven in \S 5 (Lemma 5.3). 

 In [ERT2] we showed that, under the assumptions in Theorem 1.1, 
further analysis of the heat trace associated with this problem
reduced the set of isospectral potentials for $-\Delta +q$ with 
periodic boundary
conditions to just $q(x)$ and $q(-x)$ for many choices of $q$. We plan 
to investigate the applicability of that analysis to the Dirichlet 
problem in future work.

The following result does not require the analyticity of the extended 
potential $Q$ or the irrationality of $a^2/b^2$:
\medskip
\noindent{\bf Theorem 1.3.} Suppose that $q$ and $\tilde q$ have the 
same Dirichlet spectrum on $R$, and have extensions $Q$ and $\tilde Q$ in
$C^\infty({\Bbb R}^2)$. Then, if $q(x)=q_1(x_1)+q_2(x_2)$, there are 
smooth potentials $\tilde q_1$ and $\tilde q_2$ such that $\tilde 
q(x)=\tilde q_1(x_1)+\tilde q_2(x_2)$.
\medskip
The proof of Theorem 1.3 follows the outline of the proof of the 
corresponding result, Theorem 4.1, in [ERT1]. The authors and E. Trubowitz 
proved a version of Theorem 1.3 in 1981, using the approximate 
eigenfunction construction in [ERT1, Section 3(b)], and we still plan to return 
to approximate eigenfunctions in future
work. 

As we said earlier, the generalization of Theorems 1.1 and 1.2 to 
rectangular domains in ${\Bbb R}^n$ with arbitrary assignment of
Dirichlet and Neumann boundary conditions on their faces is given in \S 
5. The constructions and arguments are sufficiently analogous to the 
two dimensional Dirichlet case that they can be presently rather 
concisely. Section 6 is devoted to
the expansion of the singularities of the wave trace at $t=|d|,\ d\in 
L\backslash 0$. This is used for the proof of Theorem 1.3 and its 
generalization to higher dimensions, Theorem 6.1. 

In \S 7 we use the Hadamard expansion of the fundamental solution to 
expand the wave
trace near $t=0$ in distributions of increasing orders of homogeneity. 
This expansion can be used to derive asymptotics for the heat trace, 
$\sum \exp(-\mu_jt)$, as $t\to 0$. While there is a huge literature on 
the asymptotics of heat traces, the relatively simple special case we 
treat
here appears to be new. The expansions reflect the singularities in the 
rectangular geometry, and contain terms which are absent for smooth 
boundaries. For instance, in dimension two with the Dirichlet condition
 $\int_Rq^2dx$ no longer appears as a spectral invariant, being 
replaced by $\int_Rq^2dx-(\pi/2)\sum q(P_i)$, where $P_i,\ i=1,..,4,$ are  the 
corners of the rectangle.

There is an extensive literature on multi-dimensional inverse spectral 
problems, and we have included a sampling in the references. In the 
setting of inverse spectral problems for $-\Delta + q$ on given domains, 
the literature includes
Guillemin-Kazdan [GuK], Guillemin [Gu1] and [Gu2], Gordon-Kappeler [GK] 
and, recently, Gordon-Sch\"uth [GS].
 
\medskip
\noindent {\bf \S 2. Fundamental Solutions, Traces and Cancellations}
\medskip
In  this section it suffices to
assume that $q$ extends to 
a sufficiently smooth
$Q\in {\Bbb R}^2$.

Let $E(t,x,y)$ be the fundamental solution for the initial value 
problem
$$u_{tt}=  \Delta u -Qu \hbox{ in } {\Bbb R}_t\times {\Bbb R}_x^2,\ 
u(0,x)=  f(x),\ u_t(0,x)=  0.\eqno{(1)}$$
It follows from the properties of the wave front sets of solutions of
the wave equation (see, for example, [H\"{o}1]) that 
$E(t,x,y)$ is an even distribution in $t$ on ${\Bbb R}$,depending 
smoothly 
on $x$ and $y$.
Here and elsewhere we follow the convention of writing distributions as 
functions, because we think it makes the manipulations easy to follow.

Using the lattice $L$ from Section 1,
the fundamental solution for (1) in $R$ with Dirichlet boundary 
conditions on
$\partial R$ can be written
$$D(t,x,y)=  $$
$$\sum_{d\in 
L}[E(t,d+x,y)-E(t,d_1-x_1,d_2+x_2,y)-E(t,d_1+x_1,d_2-x_2,y)+E(t,d-x,y)].
\eqno{(2)}$$
This construction is just the \lq\lq Method of Reflection". To check 
it, 
note first that for $(x,y)\in R\times R$ the sum is finite by finite 
speed of propagation for the wave equation. If $x_1=  0$ or $x_2=  0$ 
each term in the sum vanishes. If $x_1   =  a$, then $\pm E(t,2ma 
+a,2nb\pm 
x_2,y)$ cancels $\mp E(t,2(m+1)a-a,2nb\pm x_2,y)$, and the analogous 
cancelations occur for $x_2   =  b$. Since $E(0,x,y) =    \delta(x-y)$, 
all 
terms in the sum with $d\neq 0$ vanish for $(x,y)$ in the interior of 
$R\times R$, and we have $D(0,x,y)  =   \delta(x-y)$ and 
$\partial_t D(0,x,y)  =   0$. The symmetries of $Q$ imply that all 
terms
$E(t,d_1\pm x_1,d_2\pm x_2,y)$ are solutions of $u_{tt}-\Delta_xu 
+qu  =   0$ for $x\in R$.

Note that $D(t,x,y)$ like $E(t,x,y)$ is a distribution in  $t$ on 
${\Bbb R}$
depending smoothly on $x$ and $y$, since each of the 
terms in (2) has this property, and only finitely many of these terms 
are nonzero when $t$ lies
in a bounded interval. 
Therefore the trace $Tr( D(t))=\int_R D(t,x,x)dx$ exists
as a distribution in $t$ on ${\Bbb R}$.
Note that $\int_{-\infty}^\infty D(t,x,y)\rho(t)dt, \rho\in 
C_0^\infty$, 
is the kernel of a trace class operator 
on $L^2(R)$. Therefore
we have the 
identity
$${1\over 2}\sum [\hat \rho(\sqrt{\mu_n})+
\hat{\rho}(-\sqrt{\mu_n})] =    \int_R(\int_{\Bbb 
R}D(t,x,x)\rho(t)dt)dx$$
or more concisely in the sense of distributions on ${\Bbb R}$
$$\sum \cos(\sqrt{\mu_n}t)   =  \int_R D(t,x,x)dx  =   Tr\ 
D(t).\eqno{(3)}$$
To exploit some symmetries it is convenient to use the rectangle
$$R_0  =   \{(x_1,x_2): |x_1|\leq a,\ |x_2|\leq b\}.$$
The symmetries of $E(t,x,y)$ corresponding to the symmetries of $Q(x)$ 
imply that
$$4Tr( D(t)) = \int_{R_0}D(t,x,x)dx.\eqno{(4)}$$

>From (3) and (4) one sees that the spectrum of $H$ on $R$ determines 
$\int_{R_0}D(t,x,x)dx$. Expanding this       
using (2), we have
$$\int_{R_0}D(t,x,x)dx =     \sum_{d\in L}[\int_{R_0}E(t,d+x,x)dx 
-\int_{R_0}E(t,d_1-x_1,d_2+x_2,x)dx -$$
$$\int_{R_0}E(t,d_1+x_1,d_2-x_2,x)dx+ 
\int_{R_0}E(t,d-x,x)dx]    
=_{\hbox{def}}e^{++}(t)-e^{-+}(t)-e^{+-}(t)+e^{-
-}(t).$$
There are cancellations in $e^{-+}(t)$, $e^{+-}(t)$ and $e^{--}(t)$ 
which reduce them to sums of       
fewer terms. We will derive this reduction for $e^{--}(t)$; the 
reductions for $e^{+-}(t)$ and $e^{-+}(t)$ are
essentially corollaries of that case.

Since the periodicity of $Q$ implies $E(t,d+x,d+y) =    E(t,x,y)$ for 
all 
$d\in L$, substituting $x_1+2a$ for $x_1$ gives
$$\int_a^\infty E(t,d-x,x)dx_1  =   $$
$$\int_{-a}^\infty E(t,d-(2a,0)-x,(2a,0)+x)dx_1=     
\int_{-a}^\infty E(t,d-(4a,0)-x,x)dx_1.\eqno{(5)}$$
Likewise
$$\int_{-\infty}^{-a} E(t,d-x,x)dx_1   =  \int_{-\infty}^a 
E(t,d+(4a,0)-x,x)dx_1.\eqno{(6)}$$

Since $E(t,x,y) =    0$ when $|x-y|>|t|$,
$$\int_{-b}^b\int_{-a}^\infty E(t, (-4ma,d_2)-x,x)dx  =    0, 
\eqno{(7)}$$
when $|(-4ma,d_2)-2x|>|t|$ for $|x_2|\leq b$ and $x_1\geq -a$. Thus, 
for 
any $R$ there is an $M(R)$ such that (7)
holds for $m>M(R)$ when $|t|<R$. Thus, writing
$$\int_{-a}^a E(t,d-x,x)dx_1   =  \int_{-a}^\infty 
E(t,d-x,x)dx_1-\int_a^\infty E(t,d-x,x)dx_1,$$
(5) yields
$$\int_{-a}^a E(t,d-x,x)dx_1    = \int_{-a}^\infty 
E(t,d-x,x)dx_1-\int_{-a}^\infty E(t,d-(4a,0)-x,x)dx_1,$$

and we have the following telescoping sum formula
$$\sum_{m  =   1}^\infty\int_{-b}^b\int_{-a}^a E(t, (-2ma,d_2)-x,x)dx =    
$$
$$\int_{-b}^b\int_{-a}^\infty E(t, 
(-2a,d_2)-x,x)dx+\int_{-b}^b\int_{-a}^\infty E(t, (-4a,d_2)-x,x)dx$$
$$  =   \int_{-b}^b\int_0^\infty E(t, 
(a,d_2)-x,(a,0)+x)dx+\int_{-b}^b\int_a^\infty E(t, 
(0,d_2)-x,x)dx,\eqno{(8)}$$
where to get the last equality we replace $x_1$ by
$-a+x_1$ in the first integral and $x_1$ by $-2a+x_1$ in the second 
integral
and use the periodicity $E(t,(2a,0)-x,(2a,0)+x)=E(t,-x,x)$.

If one begins with
$$\int_{-a}^a E(t,d-x,x)dx_1  =   \int_{-\infty}^a 
E(t,d-x,x)dx_1-\int_{-\infty}^{-a} E(t,d-x,x)dx_1,$$
and uses (6), the same reasoning leads to
$$\sum_{m   =  1}^\infty\int_{-b}^b\int_{-a}^aE(t, (2ma,d_2)-x,x)dx =    
$$
$$\int_{-b}^b\int_{-\infty}^a E(t, 
(2a,d_2)-x,x)dx+\int_{-b}^b\int_{-\infty}^a E(t, (4a,d_2)-x,x)dx$$
$$   =  \int_{-b}^b\int_{-\infty}^0 E(t, 
(a,d_2)-x,(a,0)+x)dx+\int_{-b}^b\int_{-\infty}^{-a} E(t, 
(0,d_2)-x,x)dx,\eqno{(9)}$$
where to get the last equality we change $x_1$ to $a+x_1$ in the first 
integral, change $x_1$ to $x_1+2a$ in the second integral and use 
the periodicity in the second integral. 

Combining (8) and (9), we have
$$\sum_{m = -\infty}^\infty\int_{-b}^b\int_{-a}^aE(t, 
(2ma,d_2)-x,x)dx     $$
$$=\int_{-b}^b\int_{-\infty}^\infty E(t, (0,d_2)-x,x)dx+
\int_{-b}^b\int_{-\infty}^\infty E(t, (a,d_2)-x,(a,0)+x)dx
\eqno{(10)}$$

Now we can apply the same reasoning in the second variable. That gives
$$e^{--}(t)=_{\hbox{def}}\sum_{d\in L}\int_{R_0}E(t, d-x,x)dx=     
\int_{{\Bbb R}^2} E(t, -x,x)dx+\int_{{\Bbb R}^2}E(t, 
(0,b)-x,(0,b)+x)dx$$
$$+\int_{{\Bbb R}^2}E(t, (a,0)-x,(a,0)+x)dx +\int_{{\Bbb 
R}^2}E(t,(a,b)-x,(a,b)+x)dx.\eqno{(11)}$$

For $e^{-+}(t)$ and $e^{+-}(t)$ there are only cancellations in the 
sums 
in the first and second components of
the lattice vectors respectively, resulting in the formulas
$$e^{-+}(t) =    
\sum_{n =    -\infty}^\infty[ \int_{-\infty}^\infty (\int_{-b}^b E(t, 
-x_1,2nb+x_2,x)dx_2)dx_1
$$
$$+
\int_{-\infty}^\infty (\int_{-b}^b E(t, 
a-x_1,2nb+x_2,a+x_1,x_2)dx_2)dx_1],\eqno{(12)}$$
and
$$e^{+-}(t) =    
\sum_{n =    
-\infty}^\infty[\int_{-\infty}^\infty(\int_{-a}^aE(t,2na+x_1,
-x_2,x)dx_1)dx_2
$$
$$ +
\int_{-\infty}^\infty 
(\int_{-a}^aE(t,2na+x_1,b-x_2,x_1,b+x_2)dx_1)dx_2].\eqno{(13)}$$

Since $E(t,x,y)=0$ when $|x-y|>t$ the integrals over ${\Bbb R}^1$
and ${\Bbb R}^2$ in (11), (12), (13) can be written with finite limits 
of integration.
For example 
$$\int_{-b}^b\int_{-\infty}^\infty
E(t,a-x,2nb+x_2,a+x_1,x_2)dx_1dx_2 $$
$$=\int_{-b}^b\int_{4x_1^2\leq t^2-4n^2b^2}
E(t,a-x_1,2nb+x_2,a+x_1,x_2)dx_1dx_2\hbox{ and } $$
$$\int_{{\Bbb R}^2}
E(t,(a,b)-x,(a,b)+x)dx=
\int_{4|x|^2\leq t^2}E(t,(a,b)-x,(a,b)+x)dx. $$

Analogous formulas hold for other terms in (11), (12), (13). Moreover, 
for each term in $Tr(D(t))$ there is
a $t_c$ such that domain of dependence considerations, i.e. 
$E(t,x,y)=0$ for $|t|>|x-y|$, imply that
the term vanishes for $|t|<t_c$. In the next section we will see that 
these terms are real-analytic in $t$ when
$|t|>t_c$.

We can rewrite $Tr(D(t))$ in a form which has a 
geometric interpretation. We will write
$$Tr\ D(t)=\int_{R_0}D(t,x,x)dx=\sum_{m,n} D_{mn}(t),$$
where $D_{mn}(t)$ will correspond
to the lattice point $(2ma,2nb)\in L$ and have 

\noindent $t_c=2\sqrt{m^2a^2+n^2b^2}$.  When $mn\neq 0$, $D_{mn}(t)$ is 
simply
the following term from $e^{++}$
$$D_{mn}(t)=\int_{R_0} E(t,x+(2ma,2nb),x)dx,$$
Note that these terms appear in the wave trace in the case of periodic 
boundary conditions in $R_0$. 

When  $n=0$, but $m\neq 0$, we set
$$D_{m0}(t)=\int_{R_0}E(t,x+(2ma,0),x)dx$$
$$-\int_{-a}^a\int_{4|x_2|^2\leq t^2-4m^2a^2} 
E(t,x_1+2ma,-x_2,x_1,x_2)dx_1dx_2$$
$$-\int_{-a}^a\int_{4|x_2|^2\leq t^2-4m^2a^2} 
E(t,x_1+2ma,b-x_2,x_1,b+x_2)dx_1dx_2.$$
As $t\downarrow t_c$, the domain of integration in the second and third 
integrals in $D_{m0}(t)$
shrinks to $\{|x_1|\leq a,\ x_2=0\}$. Since these limiting cases 
correspond to the 
integration of $E(t,x,y)$ over $\{x=y+(2ma,0),\ |y_1|\leq a,\ y_2=0\}$
and $\{x=y+(2ma,0),\ |y_1|\leq a,\ y_2=b\}$ respectively, we will 
associate the second integral with 
the side $\{0<x_1<a,\ x_2=0\}$ of the original rectangle $R$, and the 
third with the side
$\{0<x_1<a,\ x_2=b\}$. The first integral is paired with the interior 
of $R$.  

When $m=0$, but $ n\neq 0$, we set
$$D_{0n}(t)=\int_{R_0}E(t,x+(2nb,0),x)dx$$
$$-\int_{-b}^b\int_{4x_1^2\leq t^2-4n^2b^2} 
E(t,-x_1,x_2+2nb,x_1,x_2)dx_1dx_2$$
$$-\int_{-b}^b\int_{4x_1^2\leq t^2-4n^2b^2} 
E(t,a-x_1,x_2+2nb,a+x_1,x_2)dx_1dx_2,$$
and we pair the first in integral with the interior of R, the second 
integral with 
the side $\{0<x_2<b,\ x_2=0\}$ and the third with the side $\{0<x_2<b,\ 
x_1=a\}$

Finally, when $m=n=0$ the term $D_{00}$ consists of the remaining terms 
in $Tr(D(t))$: the integral
$\int_{R_0}E(t,x,x)dx$ from $e^{++}$ paired with the interior $R$, four 
integrals from $e^{+-}$ and $e^{-+}$ with
$m=0$ and $n=0$ respectively, paired with the sides
of $R$, and the four integrals from $e^{--}$ in (11) which we pair with 
vertices
(0,0), $(a,0)$, $(0,b)$ and $(a,b)$
of $R$.  In this way we set up a one-to-one correspondence between the 
terms
in $D_{mn}(t)$ and  the k-simplices in $R$ and $\partial R$ 
corresponding to the
zero components of $(m,n)$, with $k=2$,  $k=2,1$ or $k=2,1,0$ depending 
on the number of zero 
components. This correspondence generalizes to higher dimensions (see 
Section 5).

\medskip
\noindent{\bf \S 3. Consequences of analyticity}
\medskip
In this section we draw conclusions from the analyticity of the terms 
in the sums $e^{++}(t)$, $e^{+-}(t)$, $e^{-+}(t)$ and $e^{--}(t)$. In 
[ER1]
we showed that the terms
$\int_{R_0}E(t,d+x,x)dx$
from $e^{++}(t)$ were real analytic for $|t|>|d|$ by appealing to 
results on  analytic
wave front sets. One can also prove the  results we need on the other 
terms in the trace  using analytic wave front sets. However, 
in this paper we will also use a more classical approach due to
Hadamard [H] in the form presented by  H\"ormander in [H\"o].

Let $E_+(t,x,y)$ be the forward fundamental solution for the wave
equation, i.e. the solution of $(\frac{\partial^2}{\partial t^2}-\Delta 
+Q(x))E_+= \delta_{(0,y)}$ in ${\Bbb R}_t\times {\Bbb R}^n_x$, such 
that $ E_+=0$ for $t<0$. Then $E_+$ 
 can be written as
$$E_+(t,x,y)=\sum_{\nu=0}^\infty a_\nu(x,y)e_\nu(t,|x-y|),\eqno{(14)}$$
where 
$$e_\nu(t,|x-y|)=2^{-2\nu-1}\pi^{(1-n)/2}{\Cal 
X}_+^{\nu+(1-n)/2}(t^2-|x-y|^2)$$
for $t>0, e_\nu=0$ for $t<0$.
For $a>-1$ the distribution ${\Cal X}_+^a$ is  defined  by ${\Cal 
X}_+^a(s)=(\Gamma(a+1))^{-1}s^a$ for $s>0$ and ${\Cal X}_+^a(s)=0$ for 
$s<0$. 
One defines it for all $a$ by analytic continuation in $a$; in 
particular, ${\Cal X}_+^{-k}=\delta^{(k-1)}$ for $k\in {\Bbb N}$. The
normalization of $e_\nu(t,|x-y|)$ is chosen so that  
$(\partial_t^2-\Delta) e_\nu=\nu e_{\nu-1}$ for $\nu>0$ ,  
 and
$e_0(t,|x-y|)$ is the forward fundamental solution when $Q=0$.
The coefficients $a_\nu$ are determined by the recursion
$$ \nu a_\nu +(x-y)\cdot\partial_x a_\nu +Qa_{\nu-1} -\Delta_x 
a_{\nu-1}=0 \eqno{(15)}$$
or solving (15) 
$$a_\nu(x,y)=\int_0^1s^{\nu-1}[\Delta_x 
a_{\nu-1}(y+s(x-y),y)-Q(y+s(x-y))a_{\nu-1}(y+s(x-y),y)]ds.$$
Starting with $a_0(x,y)=1$ and  using (15) repeatedly, one completes 
the construction of $E_+$. The coefficients $a_1$ and $a_2$ are given 
by
$$a_1(x,y)=-\int_0^1Q(y+s(x-y))ds\eqno{(16)}$$
and
$$a_2(x,y)=-\int_0^1s(1-s)\Delta Q(y+s(x-y))ds +{1\over 
2}(\int_0^1Q(y+s(x-y))ds)^2\eqno{(17)}$$

For general
$Q$ the sum in (14) is a singularity expansion: if one truncates the 
sum, then the difference of the fundamental solution and the truncation 
is more
regular than the last term in the truncation. However, when $Q$ is 
real-analytic one sees immediately that the $a_\nu$'s are 
real-analytic. 
Moreover, by a majorization argument one can show that the sum in (14) 
is convergent for $t$ sufficiently small as a power series in the 
variable
$(t^2-|x-y|^2)_+$. Thus we have the following convergent expansions 
for $E_+$ near $(t^2-|x-y|^2)_+=0$.
For $n$ even the function
$$A_+(t,x,y)=\sum_{\nu=0}^\infty {a_{\nu}(x,y)(t^2-|x-y|^2)^\nu\over 
2^{2\nu+1}\pi^{(n-1)/2} \Gamma(\nu +(3-n)/2)}$$
is real-analytic, and we can write
$$E_+(t,x,y)=A_+(t,x,y)(t^2-|x-y|^2)^{(1-n)/2}_+,\eqno{(18)}$$
where $(t^2-|x-y|^2)^{(1-n)/2}_+$ is interpreted in distribution sense 
as above. When $n$ is odd, one keeps the first $(n-1)/2$ terms in 
(14) and combines the remaining terms as in (18).
Let $E_-(t,x,y)=-E_+(-t,x,y)$ and let
$E_0(t,x,y)$ be the odd extension of $E_+(t,x,y)$, i.e.
$E_0(t,x,y)=E_+(t,x,y)-E_-(t,x,y)$.  Note that $E_0(t,x,y)$ is the 
solution of $(\frac{\partial^2}{\partial t^2}-\Delta +Q)u=0$
in $\Bbb R \times {\Bbb R}^n$ satisfying the initial conditions 
$E_0(0,x,y)=0,\ \frac{\partial}{\partial t}E_0(0,x,y)=\delta(x-y)$.
Therefore $E(t,x,y)=\frac{\partial}{\partial t} E_0(t,x,y)$ is
the fundamental solution for the initial value problem in (2).

Using (18) one can give convergent expansions of the trace components 
$D_{mn}(t)$ near their singularities. Consider the expression (11).  
Assuming $t>0$ and making change
of variables $x=ty$ in each integral in (11), we obtain

$$e^{--}(t)=
{\partial\over \partial t}(t\int_{|y|<1/2}{1\over 
\sqrt{1-4|y|^2}}[A_+(t,-ty,ty)+
A_+(t,(0,b)-ty,(0,b)+ty)+$$
$$+A_+(t,(a,0)-ty,(a,0)+ty)+A_+(t,(a,b)-ty,(a,b)+ty)]dy)$$
for $0<t<\epsilon$. This gives as power series expansion for 
$e^{--}(t)$ for $t>0$.

The same reasoning with the  substitution  $x_1=(t^2-4n^2b^2)^{1/2}y_1$ 
yields an expansion for 
the $n$-th term in the summation for $e^{-+}(t)$
in (12) for $t>|2nb|$, and with the 
substitution for $x_2=(t^2-4n^2a^2)^{1/2}y_2$ 
it yields an expansion for the $n$-th term in the summation for 
$e^{+-}(t)$ in (13)  for
$t>|2na|$. Since solutions of the wave equation propagate at unit 
speed, these terms are zero for $0\leq t<|2nb|$ and $0\leq t<|2na|$ 
respectively. Finally, using (18) without substitution, one gets an 
expansion for
$\int_{R_0}E(t,d+x,x)dx$ is real-analytic for $t>|d|$.

Using $d=(2ma,2nb)$, these arguments show that $ D_{mn}(t)$ is
equal to zero when $|t|<|d|$,and that it is real-analytic for
$|d|<|t|<|d|+\epsilon$. However,  we need to know that 
$D_{mn}(t)$ is real analytic on the half lines where $|t|>|d|$. For 
this we will use well-known results on analytic wave front sets. In 
particular, the following wave front calculations are based on Theorems 3.10, 
5.1 and 7.1 of [H\"o2].

Using the coordinates $(t,x,y,\tau,\xi,\eta)$ on the cotangent bundle, 
the analytic wave front set of $E$ satisfies
$$WF_A(E)\subset \{x=y + t\xi/|\xi|, \tau^2=|\xi|^2, \eta=-\xi\}.$$
Since the normal bundle to $\{(t,x+d,x)\}$ does not intersect 
$WF_A(E)$, $E(t,x,x+d)$ is well-defined as a distribution and  
$WF_A(E(t,x,x+d))\subset \{t^2=|d|^2\}$. Thus the contribution to $D_{mn}(t)$ from 
$e^{++}(t)$ is 
real analytic for $|t|\neq |d|$. 

For the terms coming from $e^{+-}(t)$ we
have to consider restrictions of $E(t,x,y)$ to 
$\{(t,-x_1,x_2+2nb,x_1,x_2)\}$ and
$\{(t,a-x_1,x_2+2nb,a+x_1,x_2)\}$. Here again the normal bundles do not 
intersect $WF_A(E)$. Thus the restrictions $E(t,-x_1,x_2+2nb,x_1,x_2)$ 
and
$E(t,a-x_1,x_2+2nb,a+x_1,x_2)$ are well-defined with analytic wave 
front sets
contained in 
$$\{ 
(\tau,\xi_1,\xi_2)=(t,x,\tau^\prime,-\xi^\prime_1+\eta^\prime_1,\xi^\prime_2+\eta^\prime_2),\tau^\prime,\xi^\prime,\eta^\prime)\in 
WF_A(E)|_{(t,x,y)}\}.$$ Since the normal bundle to a line parallel to the 
$x_1$-axis does not intersect this set when $t^2\neq 4n^2b^2$, we conclude 
that
the analytic wave front sets of 
$$\int_{\Bbb R}E(t,-s,x_2,s,x_2+2nb)ds\hbox{ and }\int_{\Bbb 
R}E(t,a-s,x_2,a+s,x_2+2nb)ds$$
are contained in $t^2=4n^2b^2$. Hence, for $n\neq 0$, $D_{0n}(t)$ is 
analytic for $|t|\neq |d|$. Analogous arguments apply to the 
contributions from $e^{-+}(t)$ and $e^{--}(t)$, and this leads to the desired 
results for $D_{m0}(t)$ and
$D_{00}(t)$.

Now we can use the analyticity results from the preceding paragraphs as 
in [ERT1] to find additional spectral
invariants. We use the irrationality of $a^2/b^2$ to conclude that 
$|d^\prime|=|d|$ for $d,d^\prime\in L$ implies
that $d^\prime=(\pm d_1,\pm d_2)$. We also assume for definiteness that 
$a<b$. Now we can recover individual terms $D_{mn}(t)$ in
 $Tr\ D(t)$ as follows.

 All  terms in $Tr\ D(t)$ except $D_{00}(t)$ are zero for $0<t<2a$ 
by domain of dependence.
So from $D(t)$ with $t$ near zero we recover $D_{00}$,  and $Tr\ D(t)$
 minus $D_{00}(t)$  is known for $t>0$ 
by
analyticity.  Continuing to reason in this way, new terms appear in the 
expansion of
$Tr\ D(t)$ at $t=|d|$ for each value of $|d|$. Since these terms are 
analytic for $t>|d|$, they are
determined by the spectrum for all $t>|d|$.
When both components of $d$ 
are nonzero, the only terms
come from $e^{++}(t)$, and we conclude that for $d_1d_2\neq 0$
$$\sum_{\pm}\int_{R_0}E(t,\pm d_1+x_1,\pm d_2+x_2,x)dx$$
is a spectral invariant. 
However, the symmetry $Q(-x_1,x_2)=Q(x_1,x_2)$ implies 

\noindent $E(t,-x_1,x_2,-y_1,y_2)=E(t,x,y)$. Thus
$$\int_{R_0}E(t,-d_1+x_1,\pm d_2+x_2,x)dx=\int_{R_0}E(t, d_1-x_1,\pm 
d_2+x_2,-x_1,x_2)dx,$$
and making the change of variables $x_1\to -x_1$, we have 
$$\int_{R_0}E(t,-d_1+x_1,\pm d_2+x_2,x)dx=\int_{R_0}E(t, d_1+x_1, \pm 
d_2+x_2,x)dx.$$
Since we also have the symmetry $Q(x_1,-x_2)=Q(x_1,x_2)$, the same 
argument applies in the variable $x_2$.
Thus the contribution from $e^{++}(t)$ simplifies to
$$ 4D_{mn}(t)=4\int_{R_0}E(t,d+x,x)dx.\eqno{(19)}$$

When $d=(d_1,0)$, $d_1\neq 0$, or $d=(0,d_2)$, $d_2\neq 0$, we get
additional contributions from $e^{+-}$ and $e^{-+}$ respectively. This 
gives the spectral invariants
$$D_{m0}(t)=\int_{R_0}E(t, d_1+x_1,x_2,x)dx_1dx_2 - 
\int_{-\infty}^\infty(\int_{-a}^aE(t,d_1+x_1,-x_2,x)dx_1)dx_2$$
$$ - 
\int_{-\infty}^\infty(\int_{-a}^aE(t,d_1+x_1,b-x_2,x_1,b+x_2)dx_1)dx_2
\eqno{(20)}$$

for $d_1=2ma\neq 0$, and
$$D_{0n}(t)=\int_{R_0}E(t,x_1,d_2+x_2,x)dx_1dx_2-\int_{-b}^b(\int_{-\infty}^\infty 
E(t, -x_1, d_2+x_2,x)dx_1)dx_2$$
$$-
\int_{-b}^b(\int_{-\infty}^\infty E(t, a-x_1, 
d_2+x_2,a+x,x_2)dx_1)dx_2\eqno{(21)}$$
for $d_2=2nb\neq 0$. In (20) and (21) we again made use of the 
observation that changing $d_i$
to $-d_i$ does not change $D_{mn}(t)$.
\medskip
\noindent {\bf \S 4. Heat Traces}
\medskip
Let $D_{mn}(t)$ be the wave trace spectral invariant from (19). Then  
$${1\over \sqrt{4\pi t}}\int_{-\infty}^\infty 
e^{-s^2/4t}D_{mn}(s)ds=\int_{R_0}
G(t,d+x,x)dx,
$$
where $G(t,x,y)$ is the fundamental solution for the initial value 
problem
$$u_t=\Delta u-Qu \hbox{ in }{\Bbb R}_+\times {\Bbb R}^2,\ 
u(0,x)=f(x).$$
Theorem 5.1 in [ERT1] gives the asymptotic behavior of $G(t,x+Nd +e,y)$ 
for $d,e\in L$ as $N\to\infty$. To describe
this we need the \lq\lq reduced potentials" 
$q_d(x)=\int_0^1Q(x+sd)ds$ and the associated fundamental solutions 
$G_d(t,x,y)$ for the initial value problems
$$u_t=\Delta u-q_du \hbox{ in }{\Bbb R}_+\times {\Bbb R}^2,\ 
u(0,x)=f(x).$$
Then 
$$|G(t,x,y)-G_d(t,x,y)|\leq {C\over 4\pi t}e^{-|x-y|^2/4t}{t\over 
|d\cdot (x-y)|}(1+{|x^\prime-y^\prime|\over t}),\eqno{(22)}$$
where $x^\prime-y^\prime$ denotes the component of $x-y$ perpendicular 
to $d$ and $C$ is
uniform in $x$ and $y$ on bounded intervals $t\in [0,T]$. In 
[ERT1] the estimate (22) is only stated for specific choices of $x$ and 
$y$, but the estimate for general $x$ and $y$
is proven. 

The structure of $G_d$ makes (22) useful. Since $Q(x+e)=Q(x)$ for $e\in 
L$, $q_d$ inherits this property, and in addition
$q_d$ is constant in the direction $d$, $q_d(x+sd)=q_d(x)$ for $s\in 
\Bbb R$. It follows that $G_d$ is the product of
a free heat kernel and a heat kernel determined by $q_d$, i.e.
$$G_d(t,x,y)={\exp[-(d\cdot (x-y))^2/4|d|^2t]\over \sqrt{4\pi 
t}}H_d(t,x,y),\eqno{(23)}$$
where $H_d(t,x+sd,y+rd)=H_d(t,x,y)$. Letting $z$ denote signed distance 
from the origin on the line $d\cdot x=0$, $H_d(t,z,w)$ is the 
fundamental solution for the initial value problem for 
$u_t=u_{zz}-q_d(z)u$. 
Letting $y=Nd+e+x$ in (22) and letting $N$ go to infinity, it follows
that for $e,d\in L$ with $d_1d_2\neq 0$ 
$\int_{R_0}G(t,x+Nd+e,x)dx$, $N\in \Bbb N$, determines 
$\int_{R_0}G_d(x+e,x)dx$ for any $d\in L\backslash 0$ and $e\in L$. 
Hence the 
invariants in (19) determine 
$$\int_{R_0}G_d(t,e+x,x)dx\eqno{(24})$$
for each $e\in L$ and $d\in L$ with $d_1d_2\neq 0$,
and the Dirichlet spectrum for $-\Delta +q(x)$ on $R$ determines
$$\sum_{e\in L} \int_{R_0} G_d(t,e+x,x)dx.\eqno{(25)}$$
This sum is the heat trace, $\sum \exp(-\lambda_n t)$, for $-\Delta 
+q_d(x)$ on $R_0$ with the
periodic boundary condition, $u(x+e)=u(x),\ e\in L$. Thus the Dirichlet 
spectrum for $-\Delta +q(x)$ on $R$ determines the periodic spectrum 
for $-\Delta +q_d(x)$ on $R_0$ for any $d\in L$ with
$d_1d_2\neq 0$.

If we choose $\delta \in L^*$ such that $\delta\cdot d=0$ and $\{\delta 
\cdot e, e\in L\}={\Bbb Z}$, then we can write the
Fourier expansion of $q_d$ in the form
$$q_d(x)=\sum_{k=-\infty}^\infty a_ke^{2\pi i k\delta\cdot x},\hbox{ 
i.e. } q_d(x)=v(\delta\cdot x)\hbox{ where }v(s)=
\sum_{k=-\infty}^\infty a_ke^{2\pi iks }.$$
In \S 3(a) of [ERT] we showed that the periodic spectrum for $-\Delta 
+q_d(x)$ on $R_0$ determines the periodic spectrum
for 
$-|\delta|^2d^2/ds^2 + v(s)$ on $[0,1]$ -- and vice versa. Note further 
that $Q(x)=Q(-x)$ and
$Q(x+d)=Q(x),\ d\in L$, imply that $q_d(x)=q_d(-x)$. Hence $a_k=a_{-k}$ 
and $v(s)=v(-s)$. The set of even potentials on $[0,1]$ with the same 
periodic spectrum is either finite or a Cantor set
(see [GT]). Thus, if $-\Delta +q$ and $-\Delta +\tilde q$ have the same 
Dirichlet spectrum on
$R$, for $d_1d_2\neq 0$  the  reduced potential $\tilde q_d$ must 
belong
to the Cantor set determined by $q_d$.

The derivation of (25) applies equally well for $d=(d_1,0),\ d_1\neq 
0$,  provided that $e_2\neq 0$, but in the remaining case, $e_2=0$, 
we need to use the invariants from (20). The same remark applies to 
$d=(0,d_2),\ d_2\neq 0$ with the invariants from (20)
replaced by those in (21). We will only give the argument for 
$d=(d_1,0)$. The
case $d=(0,d_2)$ is completely symmetric.

When $d=(d_1,0)$, $q_d(x)=q_2(x_2)$ and $G_d$ is given by
$$G_d(t,x,y)={1\over \sqrt {4\pi t}}e^{-(x_1-y_1)^2/4t}H(t,x_2,y_2),$$
where $H$ is the fundmental solution for the initial value problem
$$u_{t}=u_{x_2x_2}-q_2(x_2)u,\ u(0,x_2)=f(x_2).$$
As in (19)  the Dirichlet spectrum for $-\Delta +q$ on $R$ determines
$\int_{R_0}G_d(t,(0,e_2)+x,x)dx$ for $e_2\neq 0$, and hence it 
determines
$\int_{-b}^bH(t, x_2+e_2,x_2)dx_2$ for $e_2\neq 0$. Applying the 
asymptotic argument that lead
to (25) to (20) with $e=0$, we conclude that the Dirichlet spectrum on 
$R$ also determines
$$\int_{R_0}G_d(t,x,x)dx-\int_{-a}^a(\int_{-\infty}^\infty 
G_d(t,x_1,-x_2,x_1,x_2)dx_2)dx_1$$
$$-\int_{-a}^a(\int_{-\infty}^\infty  
G_d(t,x_1,b-x_2,x_1,b+x_2)dx_2)dx_1,$$
and hence it determines
$$\int_{-b}^bH(t,x_2,x_2)dx_2-\int_{-\infty}^\infty 
H(t,-x_2,x_2)dx_2-\int_{-\infty}^\infty H(t,b-x_2,b+x_2)dx_2.$$
Combining these two results, the Dirichlet spectrum on $R$ determines
$$[\sum_{n=-\infty}^\infty \int_{-b}^bH(t,2nb+x_2,x_2)dx_2]$$
$$-\int_{-\infty}^\infty H(t,-x_2,x_2)dx_2-\int_{-\infty}^\infty 
H(t,b-x_2,b+x_2)dx_2.\eqno{(26 )}$$
The arguments of Section 2 applied to the Dirichlet problem on $[0,b]$ 
show that
(26) is the heat trace for the Dirichlet problem for 
$-(d/dx_2)^2+q_2(x_2)$ on $[0,b]$. In
other words the Dirichlet spectrum for $-\Delta +q$ on $R$ determines 
the Dirichlet spectrum for
$-(d/dx_2)^2+q_2(x_2)$ on $[0,b]$ when $d=(d_1,0)$. Analogously, the 
Dirichlet spectrum for $-\Delta +q$ on $R$ determines the Dirichlet 
spectrum for
$-(d/dx_1)^2+q_1(x_1)$ on $[0,a]$ when $d=(0,d_2)$ 
\medskip
\noindent{\bf \S 5. More General Boundary Conditions and Higher 
Dimensions}
\medskip
The results of the preceding sections generalize to regions $R=\{x\in 
{\Bbb R}^n: 0\leq x_i\leq a_i,\ i=0,..,n\}$ with
either Dirichlet or Neumann boundary conditions on each face of 
$\partial R$. The arguments used earlier extend to
this case, because (a) the \lq\lq Method of Reflection" can be used to 
construct fundamental solutions, and (b) essential features of these 
fundamental solutions remain unchanged in the more general setting. We 
explain these points
in what follows.

In one space dimension one can check directly that the fundamental 
solutions for $u_{tt}=u_{xx}-q(x)u,\ u(0,x)=f(x),\ u_t(0,x)=0$, on the 
interval $[0,a]$ with Dirichlet or Neumann conditions on $x=0$ and 
$x=a$ can be expressed in terms of the fundamendal solution, 
$E(t,x,y)$, 
for the initial value problem $u_{tt}=u_{xx}-Q(x)u,\ u(0,x)=f(x),\ 
u_t(0,x)=0$, on the line, where
$Q$ is the extension of $q$ to ${\Bbb R}$ as before, $Q(x)=q(-x)$ for 
$-a\leq x\leq 0$ and $Q(x+2na)=Q(x)$. The formulas are 
as follows:
\medskip
\noindent for the boundary conditions $u(t,0)=u(t,a)=0$ the fundamental 
solution is
$$E_{DD}(t,x,y)= [T_{DD}E](t,x,y)=\sum_{m=-\infty}^\infty 
[E(t,2ma+x,y)-E(t,2ma-x,y)]\eqno{(27)}$$
\noindent for the boundary conditions $u(t,0)=u_x(t,a)=0$ the 
fundamental solution is
$$E_{DN}(t,x,y)= [T_{DN}E](t,x,y)=\sum_{m=-\infty}^\infty 
(-1)^m[E(t,2ma+x,y)-E(t,2ma-x,y)]\eqno{(28)}$$
\noindent for the boundary conditions $u_x(t,0)=u(t,a)=0$ the 
fundamental solution is
$$E_{ND}(t,x,y)= [T_{ND}E](t,x,y)=\sum_{m=-\infty}^\infty 
(-1)^m[E(t,2ma+x,y)+E(t,2ma-x,y)]\eqno{(29)}$$
\noindent for the boundary conditions $u_x(t,0)=u_x(t,a)=0$ the 
fundamental solution is
$$E_{NN}(t,x,y)= [T_{NN}E](t,x,y)=\sum_{m=-\infty}^\infty 
[E(t,2ma+x,y)+E(t,2ma-x,y)]\eqno{(30)}$$
These formulas can be checked as in \S 2.        
Note, for example, that in (23) we take an odd extension from
$[0,a]$ to $[-a,a]$ to satisfy the Dirichlet boundary condition at
$x=0$ and then we take an antiperiodic extension from $[-a,a]$ to
$(-\infty,+\infty)$ to satisfy the Neumann boundary condition at $x=a$.

The linear operators $T_{\alpha\beta}$ make it possible to express the 
fundamental solutions for boundary
value problems in the region $R$ above compactly. If we let 
$T_{\alpha\beta}^i$ denote $T_{\alpha\beta}$ acting
on the variable $x_i$, then the fundamental solution for $u_{tt}=\Delta 
u-qu$, $u(0,x)=f(x)$, $u_t(0,x)=0$ in $R$ with Dirichlet or Neumann
conditions on $x_i=0$ and $x_i=a_i$, $i=1,...,n$, is given by
$$E_{\alpha\beta}(t,x,y)=[T_{\alpha_1\beta_1}^1T_{\alpha_2\beta_2}^2\cdots 
T_{\alpha_n\beta_n}^n E](t,x,y),\eqno{(31)}$$
where for each $i$, $\alpha_i$ and $\beta_i$ are either $D$ or $N$. In 
(31) $E(t,x,y)$ is the fundamental solution to the initial
value problem $u_{tt}=\Delta u-Qu$, $u(0,x)=f(x)$, $u_t(0,x)=0$ in 
${\Bbb R}^n$ with $Q(\pm x_1,\pm x_2,...,\pm x_n)=q(x)$ for $x\in R$
and $Q(x_1+2m_1a_1,...,x_n+2m_na_n)=Q(x)$. In the notation 
$E_{\alpha\beta}$ in (31) we think of $\alpha$ and $\beta$ as vectors. 
Formula (31) 
reduces to (2) when
$n=2$ and one imposes the Dirichlet condition on the whole boundary. 

>From (31) one sees that the only terms in the expansion of 
$E_{\alpha\beta}(t,x,y)$ which do not have arguments of
$2m_ia_i-x_i$ for some $i$ are 
$$\sum_{(m_1,..,m_n)\in {\Bbb Z}^n}(-1)^{m_{i_1}+\cdots 
+m_{i_r}}E(t,2m_1a_1 +x_1,\dots,2m_na_n +x_n,y), \eqno{(32)}$$
where $\{i_1,\dots, i_r\}$ is the set of indices $i$ such that 
$\alpha_i\neq\beta_i$. This is the fundamental solution for
$u_{tt}=\Delta u-Qu$ in the domain
$R_0=\{x:\ -a_i\leq x_i\leq a_i\}$ with anti-periodic boundary 
conditions in the variables $x_i$ with $i=i_1,..,i_r$ and periodic
boundary conditions in the remaining variables. When we pass to traces 
as in (4), the contribution from these terms will correspond to
the wave trace for this problem. The key observation is that all other 
terms in the traces will \lq\lq telescope" as in \S 2 to 
integrals over ${\Bbb R}$ in the variables appearing in the form  
$2m_ia_i-x_i$. To see this note that the contributions of these terms 
to the 
trace 
when $\alpha_i=\beta_i$ will be summations in the index $m_i$ of the 
form
$$\cdots \sum_{m_i=-\infty}^\infty  
\int_{-a_i}^{a_i}E(t,\dots,2m_ia_i-x_i,\dots,x)dx_i\cdots,\eqno{(33)}$$
where the dots before and after the formula are intended to indicate 
summations on other indices and integrations in other variables.
The reasoning which lead to (10) reduces (33) to
$$\cdots \int_{-\infty}^\infty [E(t,\dots, -x_i,\dots,x)+E(t,\dots, 
a_i-x_i,\dots,a_i+x_i,\dots)]dx_i \cdots.\eqno{(34)}$$
When $\alpha_i\neq \beta_i$, formula (33) is replaced by 
$$\cdots \sum_{m_i=-\infty}^\infty 
(-1)^{m_i}\int_{-a_i}^{a_i}E(t,\dots,2m_ia_i-x_i,\dots,x)dx_i\cdots.\eqno{(35)}$$
To see the effect of the factor $(-1)^{m_i}$ recall the formula from 
Section 2 (rephrased in the current notation)
$$\int_{-a_i}^{a_i} E(t,...,2m_ia_i-x_i,...,x)dx_i=$$
$$\int_{-a_i}^\infty E(t,...,2m_ia_i-x_i,...,x)dx_i-\int_{-a_i}^\infty 
E(t,...,2(m_i-2)a_i-x_i,...,x)dx_i.$$
The point here is that the cancellations arise separately between the 
terms corresponding to even indices $m_i$ and those corresponding to 
odd 
indices. Hence, the terms in
(35) reduce to
$$\cdots \int_{-\infty}^\infty [E(t,\dots, -x_i,\dots,x)-E(t,\dots, 
a_i-x_i,\dots,a_i+x_i,\dots)]dx_i \cdots.\eqno{(36)}$$
\medskip
Now we can apply the arguments from Section 3 and Section 4 in this 
more general setting. Very little modification is needed. We make the 
assumptions:
\medskip
A) the extension $Q$ of $q$ to ${\Bbb R}^n$ described above is 
real-analytic on ${\Bbb R}^n$, and

B) the numbers $a_1^2,a_2^2,...,a_n^2$ are linearly independent over 
the rationals.
\medskip
\noindent We also introduce the rectangular lattice
$$L=\{(2m_1a_1,2m_2a_2,...,2m_na_n):(m_1,...,m_n)\in {\Bbb Z}^n\},$$
and its fundamental domain
$$R_0=\{x\in {\Bbb R}^n:- a_i\leq x_i<a_i,\ i=1,...,n\}.$$
Note that it follows from B) that $|d|=|d^\prime|,\ d,d^\prime\in L$ is 
equivalent to
$d_i=\pm d_i^\prime,\ i=1,...,n$. 

As was done at the end of Section 2, the distribution wave trace 
$$2^n\int_R 
E_{\alpha\beta}(t,x,x)dx=\int_{R_0}E_{\alpha\beta}(t,x,x)dx=\int_{R_0}
[T_{\alpha_1\beta_1}^1T_{\alpha_2\beta_2}^2\cdots 
T_{\alpha_n\beta_n}^nE](t,x,x)dx $$
can be represented as a sum
$$\int_{R_0}E_{\alpha\beta}(t,x,x)dx=
 \sum_{d\in L}E_{\alpha\beta}^{d}(t),\eqno{(37)}$$
where the term $E_{\alpha\beta}^{d}(t)$ corresponds to $d\in L$.  
As in Section 3 we see that each term $E_{\alpha\beta}^{d}(t)$
is real analytic when $|t|>|d|$ and equal to zero when $|t|<|d|$.
Therefore $\sum_{|\tilde{d}|=|d|}E_{\alpha\beta}^{\tilde{d}}(t)$  is
a spectral inveriant for any $d\in L$.
When $d$ has no zero components,
$$E_{\alpha\beta}^{d}(t)=\int_{R_0}E(t,d+x,x)dx.\eqno{(38)}$$
When $d_{i_1}=\cdots =d_{i_m}=0$ and $d$ has no other zero
components, the term $E_{\alpha,\beta}^{d}(t)$ is the sum of an 
integral
of form (38), $2m$  integrals associated
with the $2m$ faces $x_{i_1}=0,\ x_{i_2}=0,\cdots,\ x_{i_m}=0,
\ x_{i_1}=a_{i_1},\cdots,\ x_{i_m}=a_{i_m}$ and integrals associated
with $(n-k)$-dimensional edges, $1< k\leq m$ .
For the uniformity of notation we shall consider a face as an
$(n-1)$-dimensional edge
and each vertex of $R$ as a 0-dimensional edge.
Note that an $(n-k)$-dimensional edge is the intersection of $k$ faces.
 As in Section 2, all integrals except those of the form (38) tend to 
integrals over their
associated edges when $|t|\downarrow |d|$.

When $d$ has no 
zero components, there are the $2^n$ terms in (32) corresponding to 
$d^\prime 
$ with $|d^\prime|
=|d|$.
Since the argument leading to (19) does not depend on the number of 
variables, each of these terms (up to a {\it common} factor of $\pm 1$  
determined by the boundary conditions) equals (38).
Thus, the wave trace (37) determines the partial traces (38) for each 
$d\in L$ with no zero components.

Passing to heat traces as in Section 4, from (38) we conclude that, 
when $Nd +e$ has no zero components, the
wave trace determines
$$\int_{R_0} G(t, Nd+e+x,x)dx.$$
Taking the limit as $N\to\infty$, we obtain
$$\int_{R_0}G_d(t,x+e,x)dx,\eqno{(39)}$$
where $G_d$ is the heat kernel associated with the potential $q_d$ 
given by
$$q_d(x)=\int_0^1Q(x+sd)ds.$$
When it is possible to use this argument for all $e\in L$, we conclude 
that
the wave trace (37) determines
$$\sum_{e\in L}\int_{R_0}G_d(t,e+x,x)dx.$$
This is the trace of the fundamental solution for the heat equation on 
$R_0$
with periodic boundary conditions, i.e.  $u(t,x+d)=u(t,x),\ d\in L$, 
for the
potential $q_d$. At this point results of [ERT1] can be  applied.  
Using the notation from Theorem 1.1, the argument from pp. 668-9 of 
[ERT1] 
shows 
that for each $\delta\in S$ such that $\delta\cdot d=0$ the
invariants in (39) determine the periodic spectrum of $-\Delta + 
Q_\delta(\delta \cdot x)$.  As before the periodic spectrum of $-\Delta 
+Q_\delta(\delta\cdot x)$
determines and is determined by the periodic spectrum of 
$$-|\delta|^2{d^2\over ds^2}+Q_\delta(s)$$
on $[-1/2,1/2]$, [ERT1, Section 3(a)]. Also, $Q_\delta$ is even 
($Q_\delta(-s)=Q_\delta(s)$), and the set of even potentials with
a given periodic spectrum on $[0,1]$ is either finite or a Cantor set,
[GT]. If $\delta\in S$ has at least two nonzero components, we can 
choose
$d\in L$ with no zero components so that $\delta\cdot d=0$. Hence for 
any
$e\in L$ the vector $Nd+e$ will have no zero components for $N$ 
sufficiently
large, and we can recover the invariants in (39) for all $e\in L$. Thus
for each $\delta\in S$ with at least two nonzero components, the 
spectrum of $-\Delta +q $ on $R$ with the given boundary conditions 
determines the periodic
spectrum of $-|\delta|^2d^2/ds^2 + Q_\delta(s)$ on $[0,1]$.
\medskip
We are now left with the $\delta\in S$ with only one nonzero component, 
i.e.
$\delta =(2 a_i)^{-1}\hat e_i$ for some $i$, where $\hat e_1,...,\hat 
e_n$ is the standard basis for 
${\Bbb R}^n$. In this case $ \delta\cdot d=0$
implies $d_i=0$, and we can only recover the terms in (39) for $e\in L$ 
with
$e_i\neq 0$. To obtain the analog of the results in \S 4 in this case 
we need
to consider the sum of terms in the wave trace (37) with $t_c=|d|$ when 
the $i$-th component of $d$ is zero, but the other components are 
nonzero. This sum is given by
$$\sum_{|\tilde d|=|d|}(-1)^p[\int_{R_0}E(t,\tilde d+x,x)dx\ +$$
$$+\epsilon_1 \int_{\{|x_j|\leq a_j,\ j\neq i,\ 
|x_i|<\infty\}}[E(t,\tilde d +x-2x_i\hat e_i,x)+\epsilon_2E(t,
\tilde d + x +(a_i-2x_i)\hat e_i,x+a_i\hat{e}_i)]dx],$$
where $p$, $\epsilon_1$  and $\epsilon_2$ are determined by the 
boundary
conditions and are the same for all $\tilde d$.  In fact $\epsilon_1= 
1$ if $\alpha_i=N$ and $\epsilon_1 = -1$ if $\alpha_i=D$. Likewise 
$\epsilon_2=1$ when $\alpha_i=\beta_i$ and $\epsilon_2=-1$ when 
$\alpha_i\neq \beta_i$.
Once again the terms in this sum are independent of the sign of the 
components of $\tilde d$, and thus each term in the sum is the same. 
Passing
to the heat trace, applying the preceding with $d=Nd + e$ where 
$d_i=e_i=0$ and taking the
limit as $N\to \infty$, we conclude that the wave trace in (32) 
determines
$$\int_{R_0}G_d(t,x+e,x)dx\  +\eqno{(40)}$$
$$ +\epsilon_1 \int_{\{|x_j|\leq a_j,\ j\neq i,\ 
|x_i|<\infty\}}[G_d(t,e +x-2x_i\hat e_i,x)+\epsilon_2G_d(t,e+ x 
+(a_i-2x_i)
\hat e_i,x+a_i\hat{e}_i)]dx.$$

Now can apply a variation of the argument (from [ERT1, pp. 668-9]) used 
earlier. If $n>2$, we choose a $d^\prime\in L$ with $d_i^\prime=0$ but 
no other
nonzero components such that $d$ and $d^\prime$ are linearly 
independent. Then, taking $e=Nd^\prime+e^\prime$ with $e_i^\prime\neq 
0$ in (39) 
and
letting $N\to\infty$ we recover
$$\int_{R_0}G_{dd^\prime}(t,x+e^\prime,x)dx,\eqno{(41)}$$
where $G_{dd^\prime}$ corresponds to the potential $(q_d)_{d^\prime}$.
However, for $e^\prime$ with $e^\prime_i=0$ we make the same 
substitution in (40) and pass to the limit, getting
$$\int_{R_0}G_{dd^\prime}(t,x+e^\prime,x)dx +\eqno{(42)}$$
$$\epsilon_1 \int_{\{|x_j|\leq a_j,\ j\neq i,\ 
|x_i|<\infty\}}[G_{dd^\prime}(t,e^\prime +x-2x_i\hat 
e_i,x)+\epsilon_2G_{dd^\prime}(t,e^\prime+ 
x +(a_i-2x_i)\hat e_i,x+a_i\hat e_i)]dx.$$
We continue this until 
$(((q_d)_{d^\prime})_{d^{\prime\prime}})...=Q_\delta(\delta\cdot x)$ 
for $\delta=a_i^{-1}\hat e_i$. Letting $G_\delta$ 
denote the
heat kernel associated with the potential $Q_\delta(\delta\cdot 
x)=Q_\delta(x_i/a_i)$, we have 
$$G_\delta(t,x,y)=(4\pi 
t)^{-(n-1)/2}e^{-|x^\prime-y^\prime|^2/4t}g_i(t,x_i,y_i)$$
where $x^\prime$ and $y^\prime$ have the $i$-th components omitted, and 
$g_i$ is the fundamental solution for the initial value problem for
$u_{t}=d^2u/ dx_i^2-Q_\delta(x_i/a_i)u$. Thus the invariants in (36) 
and (37) are equivalent to
$$ \int_{-a_i}^{a_i}g_i(t,x_i+e_i,x_i)dx_i\hbox{ and }$$
$$ 
\int_{-a_i}^{a_i}g_i(t,x_i,x_i)dx_i+\epsilon_1\int_{-\infty}^\infty[g_i(t,-x_i,x_i)+\epsilon_2g_i(t,a_i-x_i,a_i+x_i)]dx_i.$$
Thus, as in the argument at the end of Section 4, undoing the 
simplification from the telescoping
terms, these invariants determine  
$$2^{-1}\int_{-a_i}^{a_i}T_{\alpha_i\beta_i}g_i(t,x_i,x_i)dx_i.\eqno{(43)}$$
Since (43) is the heat trace corresponding to the boundary condition 
$(\alpha_i,\beta_i)$ on $[0,a_i]$, we have the generalization of the 
result of Section 4:
\medskip
\noindent {\bf Theorem 5.1}. Suppose that $q$ and $\tilde q$ are 
isospectral
potentials on $R$ such that the corresponding potentials $Q$ and 
$\tilde Q$ 
on ${\Bbb R}^n$ are
real-analytic, and suppose that $\{a_1^2,\cdots,a_n^2\}$ are linearly 
independent
over $\Bbb Q$. 
Then for each $\delta \in S$ with more than one nonzero 
component  $-|\delta|^2d^2/ds^2 +Q_\delta(s)$ and $-|\delta|^2d^2/ds^2 
+\tilde Q_\delta(s)$    
have the same periodic spectrum on $[-1/2,1/2]$. For $\delta 
=(2a_i)^{-1}\hat e_i$ the operators $-d^2/dx_i^2 +Q_\delta(x_i/a_i)$ 
and 
$d^2/dx_i^2 +\tilde Q_\delta(x_i/a_i)$ have the same spectrum for the 
boundary 
condition $(\alpha_i,\beta_i)$ on $[0,a_i  ]$, $i=1,...,n$.

\noindent {\bf Remark 5.2.} When $Q_\delta=0$ for all $\delta\in S$ 
with more than one nonzero component, Theorem 5.1 gives all the 
constraints imposed by isospectrality: 
when the operators $-d^2/dx_i^2 +Q_\delta(x_i/a_i)$ and $d^2/dx_i^2 
+\tilde Q_\delta(x_i/a_i)$ have the same spectrum for the boundary 
condition 
$(\alpha_i,\beta_i)$ on $[0,a_i  ]$, $i=1,...,n$, $q$ and $\tilde q$ 
are 
isospectral. This is not what one expects when other directional 
potentials are nonzero.  In [ERT2] we used higher order terms in the 
asymptotics of $ G(t,Nd+e+x,y)$ to get additional contraints on 
isospectral 
potential and establish
rigidity modulo lattice isometries in some cases. We plan to carry this 
out for the boundary conditions considered here in the future.
\medskip
As noted in the Introduction, the set of $Q_\delta$ on $[-1/2,1/2]$ 
with a given periodic spectrum is necessarily discrete since the 
$Q_\delta$'s are even. When one works with general $L^2$ potentials, this is not 
true for the operators  $-d^2/dx_i^2 +Q_\delta(x_i/a_i)$ on $[0,a_i]$. 
See [PT, Chpt. 6] for (explicit!) examples of infinite dimensional 
manifolds of
potentials with a given Dirichlet spectrum. However, when one restricts 
the admissible potentials to those considered here, the set of 
isospectral potentials for these operators is again discrete.
\medskip
\noindent {\bf Lemma 5.3}. Assume that the potential $Q$ on ${\Bbb R}$ 
is real-analytic and satisfies $Q(x)=Q(-x)$ and $Q(x+2a)=Q(x)$. Then 
the spectrum of $d^2/dx^2 
+\tilde Q(x)$ on $[0,a]$ with either Dirichlet or Neumann boundary 
conditions at $x=0$ and $x=a$ determines the periodic spectrum of $d^2/dx^2 
+\tilde Q(x)$ on $[-a,a]$. 
\medskip
\noindent Proof. 
The reasoning given earlier, specialized to $n=1$, shows that the 
spectrum determines
$$\int_{-a}^aE(t,2ma+x,x)dx,\ m\in{\Bbb Z}\backslash 0$$
and
$$\int_{-a}^aE(t,x,x)dx+\epsilon_1\int_{-\infty}^\infty [ 
E(t,-x,x)+\epsilon_2 E(t,a-x,a+x)]dx,\eqno{(44)}$$
where $\epsilon_1$ equals 1 or -1 according to whether Neumann or 
Dirichlet conditions are imposed at $x=0$, and $\epsilon_2$ equals 1 or -1 
according to whether the boundary conditions at $x=0$ and $x=a$ are the 
same or different. Our goal is to show that as long as $Q$ is a 
real-analytic potential satisfying $Q(-x)=Q(x)$ and $Q(x+2a)=Q(x)$, the 
spectrum determines
$$\int_{-a}^aE(t,x,x)dx,\eqno{(45)}$$
and hence determines the periodic spectrum.

In one dimension the Hadamard formulas become
$$E(t,x,y)=\partial_tE_0(t,x,y),\ 
E_0(t,x,y)=E_+(t,x,y)-E_+(-t,x,y),\hbox{ and}$$
$$E_+(t,x,y)=\sum_{\nu=0}^\infty a_\nu(x,y)e_\nu(t,|x-y|),\hbox{ where 
}
e_\nu(t,|x-y|)=$$
$${1\over 2^{2\nu+1}\Gamma(\nu +1)}(t^2-|x-y|^2)^\nu_+\hbox{ for 
}t>0,\hbox{ and }E_+(t,x,y)=0\hbox{ for }t<0.$$
Hence 
$$\int_{-a}^aE(t,x,x)dx=\sum_{\nu=0}^\infty {1\over 
2^{2\nu+1}\Gamma(\nu +1)}\int_{-a}^a a_\nu(x,x)dx[\partial_t(t^{2\nu}\hbox{sgn}(t))].$$
We are going to show that the other terms in (44) contribute only even 
powers of $t$ when one expands (44) in powers of $t$, and hence one can 
recover (45) from (44).

We have (for $t>0$) 
$$\int_{-\infty}^\infty  E(t,-x,x)dx=\sum_{\nu=0}^\infty {1\over 
2^{2\nu+1}\Gamma(\nu +1)}\int_{-\infty}^\infty 
a_\nu(-x,x)[\partial_t(t^2-4x^2)^\nu_+]dx$$
$$=\sum_{\nu=1}^\infty {2\nu\over 2^{2\nu+1}\Gamma(\nu 
+1)}\int_{-1/2}^{1/2}a_\nu(ty,-ty)t^{2\nu}(1-4y^2)^{\nu-1}dy.$$
Note that in removing the term corresponding to $\nu=0$ in the last 
formula we used
$a_0(x,y)=1$. Similarly,
$$\int_{-\infty}^\infty  E(t,a-x,a+x)dx=\sum_{\nu=1}^\infty {2\nu\over 
2^{2\nu+1}\Gamma(\nu 
+1)}\int_{-1/2}^{1/2}a_\nu(a+ty,a-ty)t^{2\nu}(1-4y^2)^{\nu-1}dy$$

So to eliminate odd powers of $t$ the Maclaurin series for 
$a_\nu(-x,x)$ and $a_\nu(a-x,a+x)$ must contain only even powers of $x$. From (16) 
$$a_1(bx,x)=
-\int_0^1Q((1+(b-1)s)x)ds=\sum_{k=0}^\infty-\int_0^1{(1+(b-1)s)^k\over 
k!}Q^{(k)}(0)ds x^k$$
which contains no odd powers because $Q^{(k)}(0)=0$ for $k$ odd, and 
$$a_1(a+bx,a+x)=-\int_0^1 Q(a+(1+(b-1)s)x)ds$$
which also contains no odd powers because $Q^{(k)}(a)=0$ for $k$ odd. 
From the recursion relation (17) we have 
$$a_\nu(x,y)=\int_0^1s^{\nu-1}[\partial^2_xa_{\nu-1}(y+s(x-y),y)-Q(y+s(x-y))a_{\nu-1}(y+s(x-y)),y)]ds,$$
and 
$$a_\nu(bx,x)=$$
$$\int_0^1s^{\nu-1}[\partial^2_xa_{\nu-1}((1+(b-1)s)x,x)-Q((1+(b-1)s)x)a_{\nu-1}((1+(b-1)s)x,x)]ds,$$
Using the induction hypothesis: $a_\mu(bx,x)$  has an expansion in even 
powers of $x$ for $\mu<\nu$, we can conclude that $a_\nu(bx,x)$ has an 
expansion in even powers of $x$ (note that 
$\partial^2_xa_{\nu-1}(bx,x)=x^{-2}\partial_b^2a_{\nu-1}(bx,x)$). Likewise
$$a_\nu(a+bx,a+x)=\int_0^1s^{\nu-1}[\partial^2_xa_{\nu-1}(a+(1+(b-1)s)x,a+x)$$
$$-Q(a+(1+(b-1)s)x)a_{\nu-1}(a+(1+(b-1)s)x,a+x)]ds.$$
Hence, using the induction hypothesis: $a_\mu(a+bx,a+x)$  has an 
expansion in even powers of $x$ for $\mu<\nu$, we can conclude that 
$a_\nu(a+bx,a+x)$ has an expansion in even powers of $x$. We complete the proof 
that the Maclaurin series for $a_\nu(-x,x)$ and $a_\nu(a-x,a+x)$ must 
contain only even powers of $x$ by taking $b=-1$. Thus we have shown 
that the the spectral invariant in (44) determines the integral in (45) 
which combined with the spectral invariants preceding (44) determines the 
periodic spectrum.\qed

Lemma 5.2 combined with Theorem 5.1 gives the following generalization 
of Theorem 1.2.
\medskip
\noindent{\bf Theorem 5.4.} If $\{a_1^2,\cdots,a_n^2\}$ are linearly 
independent
over $\Bbb Q$, then for any choice $\alpha\beta$ of boundary conditions 
on the faces of $R$, any
continuous curve $q_t$ of isospectral potentials such that the 
extensions $Q_t$ are real analytic for all $t$ must be constant,
$q_t=q_0$ for all $t$.
 \bigskip
\noindent{\bf \S 6. Singularities of the wave trace and applications.}
\medskip
In this section we drop the requirements that $Q(x)$ is real-analytic 
and
$\{a_1^2,...,a_n^2\}$ is linearly independent over $\Bbb Q$, and
assume only that $Q(x)\in C^\infty({\Bbb R}^n)$.
We shall the study of the singularities of $Tr\ E_{\alpha\beta}(t)$.
We assume that $d\in L\backslash 0$ and  find the singularities of
$Tr\ E_{\alpha\beta}$ at $t=|d|$.  The case $d=0$ is discussed 
in [ER].
Using the notation of \S 5, the only terms of  $Tr\ E_{\alpha\beta}$ 
that are singular
at $t=|d|$ are

$$                                      
\sum_{|\tilde{d}|=|d|}E_{\alpha\beta}^{\tilde d}(t).        \eqno{(46)}
$$
All other terms in $Tr\ E_{\alpha\beta}$
are $C^\infty$ near $t=|d|$.
We shall find the singularity expansion of (46) at $t=|d|$, and the 
terms
of this expansion will be spectral invariants.

It follows from the Hadamard construction (see (14), (16), (17) ) that
the singularity expansion of
$$
\int_{R_0} E(t,x+d,x)dx                             \eqno{(47)}
$$
at $t=|d|$ has the form

$$
\gamma_0\int_{R_0} i dx
\frac{\partial}{\partial t} 
\chi_+^{1-n\over2}
(t^2-|d|^2)
+\gamma_1\int_{R_0}q(x)dx\frac{\partial}{\partial t} 
\chi_+^{\frac{1-n}{2}+1}
(t^2-|d|^2)
$$
$$
+\gamma_2\int_{R_0}q_d^2(x)dx\frac{\partial}{\partial t} 
\chi_+^{\frac{1-n}{2}+2}
(t^2-|d|^2)
+\cdots,                                           \eqno{(48)}
$$
where     $\gamma_i,\ 0\leq i\leq 2,$ are constants independent
of $R$ and $Q(x)$, and
$q_d(x)=\int_0^1Q(x+sd)ds$ is the reduced potential.  In deriving (48), 
we
noted that
 $Q(x)$ is periodic and smooth, and therefore
$\int_{R_0}Q(x+sd)dx=\int_{R_0}Q(x)dx=2^n\int_R q(x)dx$ and
$$
\int_{R_0}\Delta Q(x+sd)dx=0.
$$

When $d\in L$ has no zero component, $E_{\alpha\beta}^{d}(t)$
consists only of the integral (47) (up to factor $\pm 1$ 
depending on boundary conditions). Now suppose that 
$d_{i_1}=\cdots=d_{i_m}=0$ and $d$ has no other zero components.  
Then, as in Section 5, $E_{\alpha\beta}^{d}(t)$ 
in addition to (47) contains a sum of integrals associated with edges
of dimensions $n-k,\ 1\leq k\leq m$.  Consider, for example, an
integral $I_{c'}(t)$, associated with the edge $\Gamma_{c'}=\{x:  
x_{i_1}=c_{i_1},\cdot,x_{i_k}= c_{i_k}\}$, where 
the  $j$-th component of $c^\prime=(c_{i_1},...,c_{i_k})$ is either 0 
or $a_{i_j}$.
We write $x=(x',x'')$. This, as in Section 5, is a slight abuse of 
notation:
$x'=(x_{i_1},x_{i_2},\cdots,x_{i_k})$ denotes a subset of components of
$x$ - not necessarily the first $k$ components.
Near $t=|d|=|d''|$ (up to factor $\pm 1$)
$$
I_{c'}(t)=\frac{\partial}{\partial t}
\int_{\Gamma_{c'}}\int_{\Bbb R^k}E_+(t,c'-x',x'' +d'', 
c'+x',x'')dx'dx''.
                                                   \eqno{(49)}
$$
Substituting the expansion in (14) into (49), we get the expansion
$$
I_z{c'}(t)= \sum_{\nu=0}^\infty I_{c',\nu}(t),       \eqno{(50)}
$$
where 
$$
I_{c',\nu}(t)=\gamma_\nu\frac{\partial}{\partial t}
\int_{\Gamma_{c'}}\int_{\Bbb R^k}\chi_+^{\frac{1-n}{2}+\nu}(t^2-|d|^2
-4|x'|^2)a_\nu(c'-x',x''+d'',c'+x',x'')dx'dx''         \eqno{(51)}
$$
Making the change of variables 
$x'=t_dy'$ where $t_d=(t^2-|d|^2)_+^{\frac{1}{2}}$ as in Section 3, 
expanding $a_\nu(c'-t_dy',x''+d'',c'+t_dy',x'')$ in a Taylor series
in $t_dy'$, and noting that $\chi_+^{\frac{1-n}{2}+\nu}(s)$
is homogenous of degree $\frac{1-n}{2}+\nu$, we get
$$
I_{c',\nu}(t)
\approx                                                                 
\sum_{p=0}^\infty b_{c'\nu p}\frac{\partial}{\partial t}
\chi_+^{\frac{1-n}{n}+\nu+\frac{k}{2}+p}(t^2-|d|^2),           
\eqno{(52)}
$$
where
$$
b_{c'\nu 0}=\gamma_{\nu 0}(\int_{\Bbb R^k}
\chi_+^{\frac{1-n}{2}+\nu}(1-4|y'|^2)dy')
\int_{\Gamma_{c'}}a_\nu(c',x''+d'',c',x'')dx''
$$
and the
$b_{c'\nu p}$ are integrals over $\Gamma_{c'}$ of derivatives of
$a_\nu(c'-x',x''+d'',c'+x',x'')$ in $x'$ at $ x'=0$.  Using (15)
we see that $b_{c'\nu p}$ is a sum of integrals of polynomials in
$Q(x)$ and its derivatives over $\Gamma_{c'}$.
Note that integrals in $I_{c',\nu}(t)$
having odd powers of $y'$ in the Taylor formula will be zero.
 Note also that formula (52) remains true when 
$\chi_+^{\frac{1-n}{2}+\nu}$
is a distribution, i.e. when $\frac{1-n}{2}+\nu\leq -1$.  To check
this we use the representation
$$
\chi_+^{\frac{1-n}{2}+\nu}(t^2-|d|^2-4|x'|^2)
=(\frac{1}{2t}\frac{d}{dt})^N
\chi_+^{\frac{1-n}{2}+\nu+N}(t^2-|d|^2-4|x'|^2),
$$
where $\frac{1-n}{2} +\nu + N >-1$.  Now we can change variables
$x'=(t^2-|d|^2)_+^{\frac{1}{2}}y'$ as above and then apply
$(\frac{1}{2t}\frac{d}{dt})^N$ to get (52).  Collecting (6), (50), and 
(52),
we get
$$
\sum_{|\tilde{d}|=|d|}E_{\alpha\beta}^{\tilde{d}}(t)
\approx                                                 
\sum_{|\tilde{d}|=|d|}\sum_{k=1}^n\sum_{\nu,p=0}^\infty
b_{\nu k p}(\tilde{d})\frac{\partial}{\partial t}
\chi_+^{\frac{1-n}{2}+\nu+\frac{k}{2}+p}
(t^2-|d|^2)                                      \eqno{(53)}
$$
Therefore $\sum_{|\tilde{d}|=|d|}\sum_{\nu+p+\frac{k}{2}=r}
 b_{\nu k p}(\tilde{d})$
are spectral invariants for $r=0,\frac{1}{2},1,\frac{3}{2},\cdots .$

Here we will briefly describe the simplest of these spectral 
invariants.
 Since the 
coefficient $a_0$ is identically 1, $b_{\nu k p}(\tilde{d})$ with $\nu 
=0$
does not depend on $q$ and it is zero for $p>0$. Thus terms depending 
on
$q$ first appear for $r=1$. The invariant for $r=1$ is the sum of 
$b_{100}(\tilde{d})$ and $b_{020}(\tilde{d})$, and it has the form
$c_1\int_Rq(x)dx +c_2$, where $c_1$ only depends on $n$ and $c_2$ is 
determined by $R$.
Thus $\int_R q(x)dx$ is a spectral invariant as claimed in Section 1.

For $r=3/2$ the terms involving $q$ appear when $\tilde d$ has one or 
more zero components and correspond to $(\nu,k,p)=(1,1,0)$. Hence they
correspond to integrals of $q$ over the edges of $R$ of codimension 
one.
Other contributions to the invariant for $r=3/2$ 
appear only when $\tilde d$ has three or more zero components and 
correspond to $(\nu,k,p)=(0,3,0)$.

The invariants needed for the proof of Theorem 1.3 correspond to $r=2$. 
When $r=2$, contributions to the invariant come from
$(\nu,k,p)=(2,0,0),\ (1,2,0)$ and $(0,4,0)$. The $(2,0,0)$-term is the 
third term
in (48), but the $(1,2,0)$-term also depends on $q$, since it is a 
linear combination of the 
integrals of $q$ over edges of codimension 2. Since the index $k$ must 
be less than or
equal to the number of zero components of $\tilde d$, we see that
$$
\sum_{|\tilde{d}|=|d|}\int_R (q_{\tilde{d}}(x))^2 dx               
\eqno{(54)}
$$
with $d\neq 0$ will be a spectral invariant 
if $n=2$, but it will not be an invariant for $n>2$.
\medskip
\noindent {\bf Proof of Theorem 1.3.}
Suppose without loss of generality that $a\leq b$. We consider the case
$a< b$ first.  Take $v_1=(a,0)$.  Then $\tilde{d}=(\pm a,0)$ if 
$|\tilde{d}|=|v_1|$.
Therefore, using the invariance of (54), 
$\int_{R_0}q_{v_1}^2(x)dx=\int_{R_0}\tilde{q}_{v_1}^2(x)dx$.
Take $v_2=(0,b)$.  If $|\tilde{d}|=|v_2|$ then either $\tilde{d}=(0,\pm 
b)$
or $\tilde{d}=kv_1$ if $b=ka,\ k\in {\Bbb Z}$.
In both cases we have
$$
\int_{R_0}\tilde{q}_{v_2}^2(x)dx
 + C\int_{R_0}\tilde{q}_{v_1}^2(x)dx =
\int_{R_0}\tilde{q}_{v_2}^2(x)dx+C\int_{R_0}\tilde{q}_{v_1}^2dx,
$$
with $C=0$ in the first case and $C=k^2$ in the second.
Therefore 
$$
\int_{R_0}q_{v_2}^2(x)dx=\int_{R_0}\tilde{q}_{v_2}^2(x)dx.
$$
In the case when $a=b$ we have by the invariance of (54)
$$
\int_{R_0}(q_{v_1}^2+\tilde{q}_{v_2}^2)dx=
\int_{R_0}(q_{v_1}^2+q_{v_2}^2(x))dx.
$$

Take an arbitrary $d\in L$, such that $d_1 d_2\neq 0$, and consider all 
$\tilde{d}\in L$  such that $|\tilde{d}|=|d|$.  These $\tilde d$'s may 
have
the form $k_1v_1$ or $k_2v_2$ or they may have both components nonzero. 
We label the
ones with both components nonzero $\tilde d^\prime$. 
We have by (54)
$$
k_1^2\int_{R_0}\tilde{q}_{v_1}^2dx + k_2^2\int_{R_0}\tilde{q}_{v_2}^2dx
+\sum_{|\tilde{d}'|=|d|}\int_{R_0}\tilde{q}_{\tilde{d}'}^2(x)dx
$$
$$
=k_1^2\int_{R_0}q_{v_1}^2dx + k_2^2\int_{R_0}q_{v_2}^2dx
+\sum_{|\tilde{d}'|=|d|}\int_{R_0}q_{\tilde{d}'}^2(x)dx
\eqno{(55)}$$
Since  $q_{\tilde{d}'}(x)=0$ for all $\tilde{d}^\prime$, (55)
 implies $\tilde{q}_{\tilde d'}(x)=0$ for all $\tilde{d}^\prime$.  This
implies $\tilde q=\tilde{q}_1(x_1)+\tilde{q}_2(x_2)$. \ \ \ \ \ \ \ \ \ 
\ \ \ \ \ \ \ \ \ \ \ \ \ \ \
\qed

In case when the lattice $L$ satisfies the conditions of Theorem 5.1 
we can easily prove a more general result.
\medskip
\noindent {\bf Theorem 6.1.} Suppose $q(x)$ and $\tilde{q}(x)$
have the same $\alpha\beta$ spectrum on $R$ and their extensions 
$Q(x)$ and $\tilde{Q}(x)$ are sufficiently smooth in ${\Bbb R}^n$.
If $\{a_1^2,\cdots,a_n^2\}$ is linearly independent over $\Bbb Q$,
$q(x)=q_1(x_1)+\cdots +q_n(x_n)$ implies
$\tilde{q}(x)=\tilde{q}_1(x_1)+\cdots +\tilde{q}_n(x_n)$
\medskip
Proof:
\noindent
Suppose $d=(d_1,\cdots,d_n)$ has no zero components. Since we assume 
$\{a_1^2,\cdots,a_n^2\}$ is linearly independent over $\Bbb Q$,
$|\tilde{d}|=|d|$ if and only if $\tilde{d}=(\pm d_1,\pm d_2,\cdots,\pm 
d_n)$. Since none of these $\tilde d$'s have zero components,  (47),
summed over $\{\tilde d:|\tilde{d}|=|d|\}$ is a spectral invariant. 
Thus from (48) we see that (54) is a spectral invariant. 
 Since 
$q(x)=q_1(x)+\cdots+q_n(x)$ and $d$ has no zero components, we have
 $q_d(x)=0$.  Therefore by (54) $\tilde{q}_d(x)=0$ for all
$d$ with no zero components.  If $\tilde Q_\delta$ were nonzero for 
some $\delta $ with 
more than one nonzero component, we could choose a $d_0$ with no 
nonzero components such that
$\delta\cdot d_0=0$. Since this implies $\tilde q_{d_0}\neq 0$, we have 
a contradiction. Hence $\tilde Q_\delta=0$
whenever $\delta$ has more than one nonzero component, and we have 
$\tilde{q}(x)=\tilde{q}_1(x_1)+\cdots+\tilde{q}_n(x)$.
\qed

Note that in the case of periodic spectrum the analog of Theorem 6.2 
was 
proven in [ERT1] and [GK] under much weaker conditions on the lattice.
\medskip
\noindent{\bf \S 7. Heat trace asymptotics when $t\rightarrow 0$.}
\medskip

Consider the case of the Dirichlet boundary conditions in the rectangle 
$R$,
$n=2$.  For $|t|$ small we have
$$
4\ Tr\ D(t)=D_{00}(t)=
\frac{\partial}{\partial t}\int_{R_0}E_0(t,x,x)dx
$$
$$
-\frac{\partial}{\partial t}\left[\int_{-a}^a\int_{4x_2^2\leq t^2}
(E_0(t,x_1,-x_2,x_1,x_2)+E_0(t,x_1,b-x_2,x_1,b+x_2))dx_1dx_2\right.
$$
$$
\left.-\int_{-b}^b\int_{4x_1^2\leq t^2}
(E_0(t,-x_1,x_2,x_1,x_2)+E_0(t,a-x_1,x_2,a+x_1,x_2))dx_1dx_2
+e^{--}(t)\right],                                             
\eqno{(56)}
$$
where $e^{--}(t)$ has form (11).  Note that $E_0(t,x,y)$ is
an odd distribution in $t$:
$$
E_0(t,x,y)\approx \sum_{\nu=0}^\infty a_\nu(x,y)e_\nu^{(0)}(t,|x-y|),
$$
where
$$
e_\nu^{(0)}(t,|x|)=e_\nu(t,|x|)-e_\nu^-(t,|x|),
$$
$e_\nu^-(t,|x|)=0$ for $t>0,\ e_\nu^-(t,|x|)=e_\nu(-t,|x|)$ for $t< 0$.
Note that $e_\nu^{(0)}(t,|x|)$ is a distribution in $t$ depending 
smoothly
on $|x|$.  $e_\nu^{(0)}$ is odd in $t$ and homogeneous of degree 
$-1+2\nu$
when $n=2$.  In particular when $|x|=0$ we have
$$
e_0^{(0)}(t,0)=\frac{1}{2\pi}p.v.\ \frac{1}{t},\ \ 
e_1^{(0)}(t,0)=\frac{1}{4\pi}t,
\ \  e_2^{(0)}(t,0)=\frac{t^3}{24\pi}.
$$
Note that $p.v.\frac{1}{t} =\frac{d}{dt}\ln |t|$ is the only odd 
distribution
(up to a constant) homogeneous of degree -1.

Using (14), (16), (17), we recover (48) in the form
$$
\int_{R_0}E(t,x,x)dx=\frac{1}{2\pi}\frac{d^2}{dt^2}\ln |t|4ab -
\frac{1}{4\pi}\int_{R_0}Q(x)dx
+\frac{t^2}{16\pi}\int_{R_0}Q^2(x)dx+O(t^4).
$$
Consider one of four integrals in (56) associated with
sides of $R$, for example, $\frac{\partial}{\partial t}I_1(t)$,
where 
$$
I_1(t)=\int_{-b}^b\int_{4x_1^2\leq t^2} 
E_0(t,a-x_1,x_2,a+x_1,x_2)dx_1dx_2.
$$
Substituting (14), (16), (17) and making the change of variables 
$x_1=ty_1$ we get for $t>0$
$$
I_1(t)=\int_{-b}^b\int_{4y_1^2\leq 1}
\left[ \frac{1}{2\pi}\frac{1}{(1-4y_1^2)^{1/2}}
-\frac{t^2(1-4y_1^2)^{1/2}}{4\pi}\int_0^1 Q(a+ty_1-2sty_1,x_2)ds  
\right.
$$
$$
+\frac{t^4}{24\pi}(1-4y_1^2)^{3/2}
(\frac{1}{2})\left(\int_0^1Q(a+ty_1-2sty_1,x_2)ds\right)^2
$$
$$
\left.-\frac{t^4}{24\pi}(1-4y_1^2)^{3/2}
\int_0^1(1-s)s\Delta Q(a+ty_1-2sty_1,x_2)ds 
\right] dy_1dx_2 + O(t^6).
$$
We expand $Q(a+ty_1(1-2s),x_2)$ in a Taylor series in $ty_1$. Noting
that integrals containing odd powers of $y_1$ are equal to zero and
extending $I_1(t)$ for $t<0$ as an odd function we get
$$
I_1(t)=sgn\ t\left[\gamma_1 2b + t^2\gamma_2\int_{-b}^bQ(a,x_2)dx_2 
\right.
$$
$$
+t^4\left(\gamma_3\int_{-b}^bQ^2(a,x_2)dx_2 
+\gamma_4\int_{-b}^bQ_{x_1^2}(a,x_2)dx_2 \right.
$$
$$
\left.\left. +\gamma_5\int_{-b}^b\Delta Q(a,x_2)dx_2)\right) \right]
+O(t^6).
$$
Since
$Q(x_1,x_2)$ is periodic in $x_2$ we have 
$\int_{-b}^bQ_{x_2^2}(a,x_2)dx_2=0$.

Next we compute the asymptotic expansion of the four terms in 
$e^{--}(t)$
associated with vertices of $R$. Consider, for example, 
$\frac{\partial}{\partial t}I_2(t)$ where
$$
I_2(t)=\int_{4(x_1^2+x_2^2)\leq 
t^2}E_0(t,-x_1,b-x_2,x_1,b+x_2)dx_1dx_2.
$$
Making the change of variables $x_1=ty_1, x_2=ty_2$, for $t> 0$ we have
$$
I_2(t)=\int_{4|y|^2\leq 1}\left[\frac{tdy}{2\pi(1-4|y|^2)^{\frac{1}2}}
-\frac{t^3(1-4|y|^2)^{\frac{1}{2}}}{4\pi}
\int_0^1Q(ty_1-2tsy_1,b+ty_1-2tsy_2)ds\right]dy +O(t^5)
$$
$$
=\gamma_6 t +\gamma_7 t^3Q(0,b)+O(t^5).
$$
Combining the computations for all terms in (56), and evaluating the 
$\gamma_i$'s, we get
$$
\int_R D(t,x,x)dx =\frac{|R|}{2\pi}\frac{d^2}{dt^2}\ln 
|t|-\delta(t)\frac{|\partial R|}{4}+
\left(\frac {1}{4}-\frac{1}{4\pi}\int_Rq(x)dx \right)
$$
$$
+\frac{|t|}{16}\int_{\partial R}qds +t^2
\left(\frac{1}{16\pi}\int_Rq^2(x)dx-\frac{1}{32}\sum_{i=1}^4q(P_i)\right)
$$
$$
-|t|^3\left(\frac{1}{2^7}\int_{\partial R}q^2 ds -\frac{1}{2^8}
\int_{\partial R}\frac{\partial^2}{\partial n^2}q ds\right)+O(t^4).         
\eqno{(57)}
$$
Here $|R|=ab$ is the area of $R.\ |\partial R|=2a+2b$ is the length
of the boundary $\partial R,\ \frac{\partial}{\partial n}q$ is the 
normal
derivative of $q$ on $\partial R,\ P_i,\ 1\leq i\leq 4,$ are the 
vertices 
of $R$.

Denote by $G(t,x,y)$ the heat kernel in $R$ corresponding to Dirichlet 
boundary conditions. As in Section 4 we use the formula
$$
G(t,x,y)=\frac{1}{2\sqrt{\pi t}}\int_{-\infty}^\infty
D(s,x,y)e^{-\frac{s^2}{4t}}ds.                                                                 
\eqno{(58)}
$$
Applying this formula to (57) we have 
$$
\int_RG(t,x,x)dx =\frac{t^{-1}}{4\pi 
}|R|-\frac{t^{-1/2}}{8\sqrt{\pi}}|\partial R|+
\left(\frac{1}{4}-\frac{1}{4\pi}\int_Rq(x)dx\right)
+\frac{t^{1/2}}{8\sqrt{\pi}}\int_{\partial R}qds
$$
$$
+t\left(\frac {1}{8\pi}\int_Rq^2(x)dx 
-\frac{1}{16}\sum_{i=1}^4q(P_i)\right)
-t^{3/2}\left(\frac{1}{16\sqrt{\pi}}
\int_{\partial R}q^2(x)ds -\frac {1}{32\sqrt \pi}\int_{\partial 
R}\frac{\partial^2}{\partial n^2}
qds\right)+ O(t^2).                                                                 
\eqno{(59)}
$$
\bigskip
To illustrate the flexibility of the method we are using we will do one 
more example.
Consider the case of the Dirichlet boundary conditions when 
$n=3,\ R=\{0\leq x_i\leq a_i,\ i=1,2,3\}$.
Denote by $\Gamma_{c_i},i=1,2,3,$ the face $x_i=c_i$ of $R$ where 
$c_i$ is either 0 or $a_i$.  Let $\Gamma_{c_ic_j}$ be the edge,
$\Gamma_{c_ic_j}=\Gamma_{c_i}\cap\Gamma_{c_j},\ i\neq j.$ 
The trace $2^3D_{00}(t)$ is a sum of the following terms:
\medskip
\noindent $\int_{R_0}E(t,x,x)dx$ associated with the interior of $R$, 6 
integrals associated with faces $\Gamma_{c_i}$,
12 integrals associated with edges $\Gamma_{c_1c_j}$ and 
8 integrals associated with vertices $P_{c_ic_jc_k}$.

Analogously to (57) we get
$$
\int_R D(t,x,x)dx =\gamma_0\delta''(t)|R|
+\gamma_1\frac{d}{dt}p.v.\frac{1}{t}\sum_{c_i}|\Gamma_{c_i}|
+\delta(t)\left[\gamma_{21}\int_Rq(x)dx 
+\gamma_{22}\sum_{c_i,c_j}|\Gamma_{c_ic_j}|\right]
$$
$$
+\left[\gamma_{31}\sum_{c_i}\int_{\Gamma_{c_i}}qds
+\gamma_{32}\sum_{c_i,c_j,c_k}|P_{c_ic_jc_k}|\right]
+|t|\left[\gamma_{41}\int_Rq^2(x)dx 
+ \gamma_{42}\sum_{c_i,c_j}\int_{\Gamma_{c_1c_j}}qds\right]
$$
$$
+t^2\left[\gamma_{51}\sum_{\Gamma_{c_i}}\int_{\Gamma_{c_i}}q^2ds 
+ \gamma_{52}\sum_{\Gamma_{c_i}}\int_{\Gamma_{c_i}}
\frac{\partial^2}{\partial n^2}q ds
+\gamma_{53}\sum_{c_i,c_j,c_k}q(P_{c_ic_jc_k})\right]
+0(t^3).
                                                            \eqno{(60)}
$$
Here $|\Gamma_{c_i}|$ is the area of $\Gamma_{c_i},|\Gamma_{c_ic_j}|$
is the length of $\Gamma_{c_ic_j}$ and $|P_{c_ic_jc_k}|=1$ by 
definition.
Using (58) we can easily obtain the heat trace expansion from (60).
The derivation of (60) is the same as that of (57).  We only mention 
a precaution needed to handle the change of variables in
$$
\int_{\Gamma_{c_1}}\int_R e_0^{(0)}(t,2|x_1|)dx_1dx_2dx_3
$$
since $e_0^{(0)}(t,2|x_1|)$ is a distribution.
We have
$$
e_0(t,2|x_1|)=\frac{1}{2\pi}\delta(t^2-4x_1^2).
$$
  Since $e_0^{(0)}$ is 
odd in $t$, we have for $\phi\in C_0^\infty({\Bbb R}^1)$
$$
(e_0^{(0)}(t,2|x_1|),\phi)=(e_0(t,2|x|),\phi(t)-\phi(-t))=
\frac{1}{2\pi}\left(t\delta(t^2-4x_1^2),\frac{\phi(t)-\phi(-t)}{t}\right).
$$
Note that 
$\frac{\phi(t)-\phi(-t)}{t}\in C_0^\infty$ and 
$t\delta(t^2-4x_1^2)=\frac{1}{2}\frac{d}{dt}\theta(t^2-4x_1^2)$ 
where $\theta(s)=1$ when $s>0,\theta(s)=0$ when $s<0$.
Integrating $\theta(t^2-4x_1^2)$ over $x_1$ and changing variables
$x_1=ty_1$ we get
$$
\int_{{\Bbb R}}(e_0^{(0)}(t,2|x_1|),\phi)dx_1=
\frac{1}{4\pi}\left(t_+,-\frac{d}{dt} \frac{\phi(t)-\phi(-t)}{t}\right)
=\frac{1}{4\pi}\left(\theta(t),\frac{\phi(t)-\phi(-t)}{t}\right)=
\frac{1}{4\pi}(p.v.\frac{1}{t},\phi).
$$
Similar arguments   show that $(e_0^{(0)}(t,0),\varphi(t))=
\frac{1}{2\pi}\phi'(0)=-\frac{1}{2\pi}(\delta',\phi)$

With a little more effort one can compute the coefficients in (60). 
They are
$\gamma_0=-(2\pi)^{-1}$, $ \gamma_1=-(8\pi)^{-1}$, 
$\gamma_{21}=-(4\pi)^{-1}$, $
\gamma_{22}=(16)^{-1}$, $ \gamma_{31}=(16\pi)^{-1}$, $ 
\gamma_{32}=-(64)^{-1}$, 
$\gamma_{41}=(32\pi)^{-1}$, 
$\gamma_{51}=-(64)^{-1}$, $\gamma_{52}=(128\pi)^{-1}$ and $
\gamma_{53}=(128)^{-1}$.
Applying (58) to the expansion in (60), we get the expansion of the
heat trace
$$\sum e^{-\mu_jt}=$$
$${1\over 8(\pi t)^{3/2}}|R| -{1\over 16\pi
t}\sum_{c_i}|\Gamma_{c_i}| +{1\over t^{1/2}}\left[{-1\over
8\pi^{3/2}}\int_Rq(x)dx +{1\over
32\pi^{1/2}}\sum_{c_i,c_j}|\Gamma_{c_ic_j}|\right]
$$
$$
+\left[{1\over 16\pi}\sum_{c_i}\int_{\Gamma_{c_i}}qds
-{1\over 64}\sum_{c_i,c_j,c_k}|P_{c_ic_jc_k}|\right]
+t^{1/2}\left[{1\over 16\pi^{3/2}}\int_Rq^2(x)dx
- {1\over 32\pi^{1/2}}\sum_{c_i,c_j}\int_{\Gamma_{c_1c_j}}qds\right]
$$
$$
+t\left[-{1\over 32\pi}\sum_{c_i}\int_{\Gamma_{c_i}}q^2d\sigma
+ {1\over 64\pi}\sum_{c_i}\int_{\Gamma_{c_i}}
\frac{\partial^2}{\partial n^2}q d\sigma
+{1\over 64}\sum_{c_i,c_j,c_k}q(P_{c_ic_jc_k})\right]
+0(t^{3/2}).$$
\medskip
These computations can, of course, be carried out for arbitrary
assignments of Dirichlet and Neumann conditions on the sides 
$\Gamma_{c_i}$.
The only changes are as follows. The integrals in $D_{00}$ 
corresponding to $\Gamma_{c_i}$ with
the Dirichlet condition are preceded by minus signs while those 
corresponding to $\Gamma_{c_i}$ are
preceded by plus signs. The integrals corresponding to edges 
$\Gamma_{c_ic_j}$ with the same
boundary condition on $\Gamma_{c_i}$ and $\Gamma_{c_j}$ are preceded by 
plus signs, and the 
others have minus sign. Finally the integrals corresponding to vertices 
are preceded by $(-1)^m$ where
$m$ is the number of adjoining faces with the Dirichlet condition.

\medskip
\centerline{\bf References}
\medskip
\noindent [CH] Courant, R. and Hilbert, D., {\it Methods of 
Mathematical Physics, Vol. II}, Interscience, New York, 1962.
\medskip
\noindent [ERT1] Eskin, G., Ralston, J. and Trubowitz, E., On 
isospectral periodic potentials in ${\Bbb R}^n$, Comm. Pure and
App. Math. {\bf 37}(1984), 647-676.
\medskip 
\noindent [ERT2 ]Eskin, G., Ralston, J. and Trubowitz, E., On 
isospectral periodic potentials in ${\Bbb R}^n$,II, Comm. Pure and
App. Math. {\bf 37}(1984), 715-753.
\medskip
\noindent [GT] Garnett, J. and Trubowitz, E., Gaps and bands of 
one-dimensional periodic Schr\"od-

\noindent inger operators, Comm. Math. Helvetici
{\bf 59}(1984), 258-312.
\medskip
\noindent [G] Gordon, C., Survey of isospectral manifolds, {\it 
Handbook of Differential Geometry}, Vol. 1, 747-778, North Holland, 
2000.
\medskip 
\noindent [GK] Gordon, C. and Kappeler, T., On isospectral potentials 
on tori, Duke Math. J. {\bf 63}(1991), 217-233.
\medskip
\noindent [GS] Gordon, C. and Sch\"uth, D. Isospectral potentials and 
conformally equivalent isospectral metrics on spheres, balls and Lie
groups, J. Geometric Anal. {\bf 13}(2003), 300-328.
\medskip
\noindent [Gu1] Guillemin, V., Spectral theory on $S^2$: some open 
questions, Adv. Math. {\bf 42}(1981), 283-298.
\medskip
\noindent [Gu2] Guillemin, V., Inverse spectral results on 
two-dimensional tori, JAMS {\bf 3}
(1990), 375-387.
\medskip
\noindent [GuK] Guillemin, V. and Kazdhan, D., Some inverse spectral 
results for negatively curved 2-manifolds, Topology {\bf 19}(1980), 
301-312.
\medskip
\noindent [H] Hadamard, J. {\it Le probl\`eme de Cauchy et les 
\'equations aux d\'eriv\'ees partielles lin\'eaires hyperboliques}, 
Hermann, 
Paris, 1932.
\medskip
\noindent [Ha] Hald, O. and McLaughlin, J., Inverse nodal problems: 
finding the potential from nodal lines, Mem. Amer. Math. Soc. {\bf 
119}(1996), no. 572.
\medskip
\noindent [H\"o1] H\"ormander, L., {\it The Analysis of Linear Partial 
Differential Equations, III}, Springer-Verlag, Vienna, 1985.
\medskip
\noindent [H\"o2] H\"ormander, L., Uniqueness theorems and wave front 
sets for solutions of linear differential equations with analytic 
coefficients, Comm. Pure Appl. Math {\bf 24}(1971), 671-704.
\medskip
\noindent [KT] Kn\"orrer, H. and Trubowitz, E., A directional 
compactification of the complex Bloch variety, Comm. Math. Helvetici 
{\bf 
65}(1990), 114-149.
\medskip
\noindent [Mc] McLaughlin, J., Solving inverse problems with spectral 
data, {\it Surveys of Solution Methods for Inverse Problems}, 169-194, 
Springer-Verlag, Vienna, 2000.
\medskip
\noindent [PT] Popov, G. and Topalev, P., Liouville billiard tables and 
an inverse spectral result, Ergodic Thy. and Dyn. Sys. {\bf 23}(2003), 
225-248.
\medskip
\noindent [Z] Zeldich, S., Spectral determination of analytic, 
bi-axisymmetirc plane domains, Geometric Funct. Anal. {\bf 10}(2000), 
628-677.

\end